\documentclass[12pt, reqno]{amsart}
\usepackage{graphicx}
\usepackage{float}
\usepackage{epstopdf}
\usepackage{tabularx}
\usepackage{caption}
\usepackage{times}
\usepackage{amsmath}
\usepackage{amsfonts}
\usepackage{amssymb}
\usepackage{multirow}
\usepackage{xmpmulti}
\usepackage{booktabs}
\usepackage[colorlinks, linkcolor=red, anchorcolor=blue, citecolor=blue]{hyperref}
\usepackage{fancyhdr}
\newtheorem{theorem}{Theorem}

\newtheorem{example}{Example}
\newtheorem{remark}[theorem]{Remark}

\def\D{\bf D}
\def\C{\bf C}
\makeatletter
\@namedef{subjclassname@1991}{2020 Mathematics Subject Classification}

\begin{document}

\title{AFD Types Sparse Representations vs. \\the Karhunen-Lo\`eve Expansion for Decomposing Stochastic Processes}

\author{Tao Qian}
\address[Qian]{Macau Center for Mathematical Sciences, Macau University of Science and Technology, Macau.}
\email{tqian$\symbol{64}$must.edu.mo}
\author{Ying Zhang}
\address[Zhang]{Macau Center for Mathematical Sciences, Macau University of Science and Technology, Macau.}
\email{cnuzhangying$\symbol{64}$163.com}
\author{WanQuan Liu}
\address[Liu]{School of Intelligent Systems Engineering, Sun Yat-sen University, Shenzhen, China.}
\email{liuwq63$\symbol{64}$mail.sysu.edu.cn}
\author{Wei Qu*}
\address[Qu]{College of Sciences, China Jiliang University, China.}
\email{quwei2math$\symbol{64}$qq.com}

\thanks{*Corresponding author. The study is supported by the Science and Technology Development Fund, Macau SAR (File no. 0123/2018/A3). }

\begin{abstract}
This article introduces adaptive Fourier decomposition (AFD) type methods, emphasizing on those that can be applied to stochastic processes and random fields, mainly including stochastic adaptive Fourier decomposition and stochastic pre-orthogonal adaptive Fourier decomposition. We establish their algorithms based on the covariant function and prove that they enjoy the same convergence rate as the Karhunen-Lo\`eve (KL)  decomposition. The AFD type methods are compared with the KL decomposition. In contrast with the latter, the AFD type methods do not need to compute eigenvalues and eigenfunctions of the kernel-integral operator induced by the covariance function, and thus considerably reduce the computation complexity and computer consumes. Various kinds of dictionaries offer AFD flexibility to solve problems of a great variety, including different types of deterministic and stochastic equations. The conducted experiments show, besides the numerical convenience and fast convergence, that the AFD type decompositions outperform the KL type in describing local details, in spite of the proven global optimality of the latter.
\end{abstract}



\maketitle
{\center
\section*{Declaration of interests}
The authors declare that they have no known competing financial interests or personal relationships that could have appeared to influence the work reported in this paper.}
\section{Introduction}
For the self-containing purpose this section will introduce adaptive Fourier decomposition (AFD) type sparse representations with emphasis on stochastic AFDs (\cite{QQ,QianSAFD}).
We are based on a dictionary $\mathcal{D}$ of a complex Hilbert space $\mathcal{H}.$  By definition, a dictionary of $\mathcal{H}$ consists of a class of unimodular elements whose linear span is dense in $\mathcal{H}.$ The formulation we adopt is that $\mathcal{H}$ is the $L^2$-space of complex-valued functions on a manifold $\partial {\D},$ being the boundary of ${\D},$ where $\D$ itself is an open and connected domain, called a region, in an Euclidean space. A process in a finite time interval may be well treated as defined on the unit circle $\partial {\D}=\{e^{it}\ |\ t\in [0,2\pi)\}.$ If a process is defined in the entire time range, we use the format that $\partial \D$ be the whole real line and $\D$ itself the upper-half plane. For general random fields we may fit the problem to the cases where $\D$ being the solid ball in $\mathbb{R}^n,$ denoted as ${\bf D}_1,$ or the upper-half spaces $\mathbb{R}^{n+1}_+$, denoted as ${\bf D}_2,$ etc. The concerned contexts may be fit into the following two structures, the coarse structure and the fine structure, as described below.\\

\noindent (i). The coarse structure: $\mathcal{H}=L^2(\partial \D)$ and the elements of the dictionary $\mathcal{D}$ are indexed by all $q$ in $\D.$ Examples for this model include $\mathcal{D}$ being the collection of the Poisson kernels in the unit disc or the unit ball ${\D}_1$, or those in the upper-half space ${\D}_2.$  Alternatively, $\mathcal{D}$ can be the collection of the heat kernels in the upper-half space ${\D}_2.$ The upper-half space can also be associated with general dilated and translated convolution kernels (\cite{Qu-Sps-represt-dirac}).  \\

\noindent (ii). The fine structure: Certain functions defined on a region $\D$ may constitute a reproducing kernel Hilbert space (RKHS), be denoted by $H_K,$ or $H^2(\D),$ and called the Hardy space of the context, where $K: \D\times\D\to \C$ is the reproducing kernel, $K_q(p)=K(p,q), p,q\in \D.$ In such case the boundary values, or the non-tangential boundary limits, as often occur in harmonic analysis, span a dense class of $\mathcal{H}=L^2(\partial \D)$ (as a new notation),  or sometimes just dense in a proper subspace of $\mathcal{H}$. In the case the boundary value mapping induces an isometry, or a bounded linear operator, between $H^2(\D)$ and $H^2(\partial \D).$ The related theory and examples are contained in \cite{Saitoh_RKHS-2016}, \cite{Qian2020HHK}, \cite{Qu-Sps-represt-dirac}, as well as in \cite{YangF}. The relevant literature address various types of RKHSs, or Hardy spaces in such setting, including, for instance, the $H^2$ space of complex holomorphic functions in the unit disc ${\D}_1,$ and the $h^2$ space of harmonic functions in the upper-half Euclidean space, precisely,
\[H^2({\D}_1)=\{f:{\D}_1\to {\C}\ |\ f(z)=\sum_{k=0}^\infty c_kz^k, \|f\|^2=\sum_{k=0}^\infty |c_k|^2<\infty\};\]
the harmonic $h^2$-spaces in the upper-half space ${\D}_2,$
\[ h^2({\D}_2)=\{u:{\D}_2\to {\C}\ |\ \Delta u=0\ {\rm on\ }
{\D}_2, \sup_{(t,y)\in \Gamma^\alpha_x} |u(t,y)|\in L^2(\mathbb{R}^n)\},\]
where $\Gamma^\alpha_x$ is the orthogonal $\alpha$-cone in $\mathbb{R}_{+}^{1+n}$ with tip at $x\in \mathbb{R}^n;$ and similarly the harmonic $h^2$-space of the unit $n$-ball in $\mathbb{R}^{1+n};$ and the heat kernel Hardy space $H^2_{heat}(\mathbb{R}_{+}^{1+n}).$ In the present paper all such RKHSs are denoted $H^2(\D).$\\

As a property of reproducing kernel, in each of the fine cases, the parameterized reproducing kernels $K_q, q\in {\D},$ form a dictionary. To simplify the terminology we, in this paper, also call a dictionary element as a kernel.
Based on the above formulation we now briefly review a few sparse representation models belonging to the AFD type. Unless otherwise specified, the norm $\|\cdot\|$ and the inner product notation $\langle \cdot,\cdot \rangle$ always mean those of the underlying complex Hilbert space $\mathcal{H}.$ \\

\noindent (a).  AFD, or Core AFD:  In the $\mathcal{H}=H^2({\D}_1)$ context we have, being based on the Szeg\"o kernel dictionary,
\[ \mathcal{D}_1\triangleq\{k_a(z)\}_{a\in {\D}_1}\triangleq\{\frac{1}{1-\overline{a}z}\}_{a\in {\D}_1}.\] For any given $f\in H^2({\D}_1)$ we have a greedy selections of the parameters
\begin{eqnarray}\label{MSPD3} a_k=\arg \max \{ |\langle f_k,e_a\rangle |\ |\ a\in {\D}_1\},\end{eqnarray}
where  $f_k$ are the reduced remainders with $f_1=f,$ obtained through inductively  using the \emph{generalized backward shift operator}:
\[ f_k(z)=\frac{f_{k-1}(z)-\langle f_{k-1},e_{a_{k-1}}\rangle e_{a_{k-1}}(z)}{\frac{z-a_{k-1}}{1-\overline{a}_{k-1}z}},\]
where $e_a=k_a/\|k_a\|, \|k_a\|=1/\sqrt{1-|a|^2},$ is the normalized Szeg\"o kernel, and the validity of the Maximal Selection Principle (\ref{MSPD3}) (MSP) is a consequence of the Boundary Vanishing Condition (BVC)
\[ \lim_{|a|\to 1}\langle f_k,e_a\rangle =0.\]
Then there follows (see \cite{QWa})
\[f(z)=\sum_{k=1}^\infty \langle f_k,e_{a_k}\rangle B_k(z),\]
where $\{B_k\}$ is automatically an orthonormal system, called the Takenaka-Malmquist (TM) system, given by
\begin{eqnarray}\label{TM} B_k(z)=e_{a_k}(z)\prod_{l=1}^{k-1}\frac{z-a_{l}}{1-\overline{a}_{l}z}.\end{eqnarray}
Whether the associated TM system is a basis depends on whether $\sum_{k=1}^\infty
(1-|a_k|)=\infty.$ In both the basis and the non-basis cases the above infinite series converges in the $L^2(\partial {\D}_1)$ norm sense on the boundary with the convergence rate
\begin{eqnarray}\label{rate} \|f-\sum_{k=1}^n \langle f_k,e_{a_k}\rangle B_k\|\leq \frac{M}{\sqrt{n}},\end{eqnarray}
where
\[ M=\inf\{\sum_{k=1}^\infty |c_k|\ |\ f(z)=\sum_{k=1}^\infty c_ke_{b_k}, b_k\in {\D}_1\}\]
(see \cite{QWa} or \cite{Te2}).
Being compared with the usual greedy algorithm, Core AFD was originated by the purpose of finding meaningful positive frequency decomposition of a signal. It was found that analytic phase derivatives of inner functions provide the source of such meaningful positive frequencies (\cite{QCL,Qianinner}). Apart from this background, AFD addresses attainability of globally maximal energy matching pursuit. The attainability in various contexts of kernel approximation further motivated the multiple kernel concept (see below).
For Core AFD the reader is referred to \cite{QWa,Q2D}. For general greedy algorithms and convergence rate results the reader is referred to \cite{DT} and \cite{Te2}. Core AFD has the following variations and generalizations.\\

\noindent (b) Unwinding AFD: The inductive steps in (a) give rise to the relation
\[ f(z)=\sum_{k=1}^{n-1} \langle f_k,e_{a_k}\rangle B_k+g_{n}(z),\]
where
\[ g_n(z)=f_n(z)\prod_{l=1}^{n-1}\frac{z-a_{l}}{1-\overline{a}_{l}z}\]
is the \emph{standard $n$-orthogonal remainder}.
If, before doing the optimal energy matching pursuit for $f_n,$ one first gets out the inner function factor of $g_n,$ and then perform a maximal parameter selection, one obtains what is called unwinding AFD (\cite{Qian2010}), which converges considerably faster than Core AFD and, again, with positive frequencies. Coifman et al. studied a particular type of unwinding AFD that was named as Blaschke unwinding expansion (\cite{Nahon, Coifman1, CP}).\\

\noindent (c) Based on generalizations of Blaschke products the AFD theory is extended to several complex variables contexts, and to matrix-valued functions (see \cite{ACQS1,ACQS2}). \\

\noindent (d) In a general context, however, there may not exist a Blaschke product theory, and the Gram-Schmidt (GS) orthogonalization of the parameterized reproducing kernels or the dictionary elements may not have explicit and practical formulas as those in the TM system (\ref{TM}). Some relations held for the Core AFD case, however, motivate what we call as pre-orthogonal AFD (POAFD), the latter is equivalent with Core AFD when the context reduces to the classical Hardy space. Precisely speaking, due to orthogonality of the TM system, there hold
\[ \langle f_n,e_{a_n}\rangle=\langle g_n,B_n\rangle=\langle f,B_n\rangle.\]
The above relations hint the notion $B_n^a$ within the $n$-tuple
$(B_1,\cdots,B_{n-1},B_n^a),$ where the latter being the GS orthonormalization of $(B_1,\cdots,B_{n-1},k_a)$ (\cite{QianSAFD}), and $a$ is left to be determined.

Below for a general Hilbert space $\mathcal{H}$ we use the notations ${\D}, q, K_q, E_q, E_n$ to replace ${\bf D}_1, a, k_a, e_a, B_n$ specially used for the Core AFD case.
We now introduce Pre-orthogonal Maximal Selection Principle (\ref{Pre-MSP}) (Pre-MSP) by which Core AFD is generalized to any Hilbert space with a dictionary satisfying BVC (or called a \emph{BVC dictionary}):
Select
\begin{eqnarray}\label{Pre-MSP}
q_n=\arg\sup \{|\langle f,E_n^q\rangle|\ |\ q\in {\D} \},
\end{eqnarray}
where $\{E_1,\cdots,E_{n-1},E_n^q\}$ is the GS orthogonalization of
$\{E_1,\cdots,E_{n-1},\tilde{K}_{q}\},$
where for each $j: 1\leq j\leq n,$
 \[\tilde{K}_{q_j}=\left[\left(\frac{\partial}{\partial \overline{q}}\right)^{(l(j)-1)}K_q\right]_{q=q_j},\ j=1,2,\cdots,n,\]
 where $\frac{\partial}{\partial \overline{q}}$ is a directional derivative with respect to $\overline{q},$ and $l(j)$ is the multiple of $q_j$ in the $j$-tuple $(q_1,\cdots,q_j), 1\leq j\leq n.$ $\tilde{K}_{q_j}$ is sometimes briefly denoted as
 $\tilde{K}_{j}$ and called the \emph{multiple kernel} with respect to $q_j$ in the $j$-tuple $(q_1,\cdots,q_j).$ With this notation we have $E_q=E^q_1,$ as used in (\ref{below}). In order to attain the supreme value at each of the parameter selections, one must allow repeating selections of parameters, corresponding to multiple kernels (\cite{Qian2018,CQT}). Recall that the basic functions in a TM system correspond to the GS orthogonalizations of the involved multiple Szeg\"o kernels (\cite{QianSAFD}). In the classical Hardy space case POAFD reduces to AFD. We note that the latter is the only case in which the GS process generates out such nice and practical formulas (see (\ref{TM})).
 To perform GS orthogonalization is to compute
   \begin{eqnarray}\label{GSmultiple}
 E_{n}
 =\frac{\tilde{K}_{q_{n}}-\sum_{k=1}^{n-1}\langle \tilde{K}_{q_{n}},E_k\rangle E_k}
 {\sqrt{\|\tilde{K}_{q_n}\|^2-\sum_{k=1}^{n-1}|\langle \tilde{K}_{q_{n}},E_k\rangle|^2}}.\end{eqnarray}
   With such formulation Core AFD can be extended to various contexts in which a practical Blaschke product theory may not be known or may not exist.
   Significant generalizations include POAFD for product dictionary (\cite{Q2D}), POAFD for quaternionic space (\cite{QSW}), POAFD for multivariate real variables in the Clifford algebra setting (\cite{WjQ}), POAFD for weighted Bergman and weighted Hardy spaces (\cite{qu2018,qu2019}), and most recently sparse representations for the Dirac $\delta$ function \cite{Qu-Sps-represt-dirac}, and etc. We note that the maximal selection principle of POAFD guarantees that, in the optimal one-step selection category, the greediest algorithm, and, in particular, more greedy than what is called \emph{orthogonal greedy algorithm } (\cite{Te2}) due to the relation: For $q\ne q_,\cdots,q_{n-1},$
   \begin{eqnarray}\label{below} |\langle g_n,E^{q}_n\rangle|= \frac{1}{\sqrt{1-\sum_{k=1}^{n-1}|\langle E_q,E_k\rangle|^2}}|\langle g_n,E_1^q\rangle|\ge |\langle g_n,E_q\rangle|,\end{eqnarray}
   the ending term being what is used for a maximal parameter selection of the orthogonal greedy algorithm.\\

 \noindent (e) $n$-Best AFD: A manipulation of a single optimal parameter selection with Core AFD or POAFD in (\ref{MSPD3}) or (\ref{Pre-MSP}), respectively, is $n$-best AFD, also called $n$-\emph{best kernel approximation}, formulated as finding $(q_1,\cdots,q_n)$ such that
 \begin{eqnarray}\label{nbest}
 & &\|P_{{\rm span}\{\tilde{K}_{q_1},\cdots,\tilde{K}_{q_n}\}}(f)\|\nonumber \\
 &=&\sup\{\|P_{{\rm span}\{\tilde{K}_{p_1},\cdots,\tilde{K}_{p_n}\}}(f)\|\ |\ p_1,\cdots, p_n \in {\D} \},
 \end{eqnarray}
 where we use $P_X(f)$ as the projection of $f$ into a linear subspace $X.$ In the classical Hardy space case ${\D}={\D}_1$ the problem is equivalent with the one of finding best approximation by rational functions of degree not exceeding $n.$ The existence part of the rational approximation problem has long been solved (see \cite{WQ2020} and the references therein), the algorithm, however, is left open until now (\cite{Qiancyclic}). In general, for a Hilbert space with a BVC dictionary asserting existence of a solution for the $n$-best problem is by no means easy. The existence result for the classical complex Hardy space case has recently been re-proved by using a new approach based on the maximal module principle of complex analytic functions (\cite{WQ2020}). This progress allows to generalize the existence of an $n$-best solution to weighted Bergman spaces (\cite{QQ}), and further to a wide class of RKHSs (\cite{Qiannbest}): Under a set of commonly used conditions existence of an $n$-best approximation is proved for a large class of RKHSs consisting of certain analytic functions in the unit disc.  In the upper-half of the complex plane there is a parallel theory.

 \noindent (f) The most up-date developments of AFD is stochastic AFD (SAFD) and stochastic pre-orthogonal AFD (SPOAFD) (\cite{QianSAFD}). The former is precisely for the classical complex Hardy space context which, as mentioned, with the convenience of the TM system, and the latter is for general cases with a BVC dictionary, including a wide class of functional spaces. The purpose of the present paper is to introduce SAFD and SPOAFD to the study of, and practice with stochastic processes and random fields.
 \newtheorem{Def}{Definition}
 \begin{Def} (\cite{LPS,QianSAFD})
	Suppose ${\D} \subset \mathbb R ^n.$ The Bochner type space $ L^2(\Omega, L^2({\D})) $ is defined to be the Hilbert space consisting of all the $ L^2({\D}) $-valued random variables $f: {\D} \times \Omega \rightarrow \mathbb{C}$ that satisfy
	\begin{eqnarray}
	\|f\|_{L^2(\Omega, L^2({\D}))}^2&\triangleq&
	\int_{\Omega} \int_{{\D}} |f(x, \omega)|^2 dx d\mathbb{P}\nonumber \\
&=& E_\omega\|f_\omega\|^2_{L^2({\D})}<\infty,
	\end{eqnarray}
where $f_\omega(x)=f(x,\omega).$
\end{Def}
\def\N{\mathcal N}
\noindent
For brevity, we also write $ \mathcal{N}=L^2(\Omega,L^2({\D})).$ The SAFD theory developed in \cite{QianSAFD} is for the same space in the unit disc ${\D}_1$ but defined in terms of the Fourier expansion with random coefficients, being equivalent with the above defined due to the Plancherel Theorem. Besides SAFD, \cite{QianSAFD} also develops SPOAFD for general stochastic Hilbert spaces with a BVC dictionary.

 SAFD (identical with SAFD2 in the terminology of \cite{QianSAFD}) concerning the complex Hardy space $H^2({\D}_1)$ with the Szeg\"o kernel dictionary $\mathcal{D}_1,$ precisely corresponds to $ \mathcal{N}=L^2(\Omega,H^2({\D}_1)).$
 For $f\in \N,$ there holds, for a.s. $\omega\in \Omega,$ $f_\omega\in H^2({\D}_1),$ and the stochastic MSP (SMSP) is, for $f_{\omega}(e^{it})=f(e^{it},\omega),$
\begin{eqnarray}\label{MSPD1} a_k=\arg \max \{\mathbb{E}_{\omega} |\langle f_{\omega},B_k^a\rangle |^2\ |\ a\in {\D}_1\},\end{eqnarray}
where for the previously known $a_1,\cdots,a_{k-1},$
\begin{eqnarray}\label{no} B_k^a(z)=e_{a}(z)\prod_{l=1}^{k-1}\frac{z-a_{l}}{1-\overline{a}_{l}z}\end{eqnarray} is the $k$-th term of the TM system with the parameter $a\in {\D}_1$ to be optimally determined.
 Existence of an optimal $a_k$ is proved in \cite{QianSAFD}. Then for $B_k=B_k^{a_k}$ under such an optimal $a_k$ the consecutively determined TM system gives rise to an expansion of $f$
 \[ f(e^{it},\omega)\stackrel{\mathcal{N}}{=}\sum_{k=1}^\infty \langle f_\omega,B_k\rangle B_k(e^{it})\]
 for a.s. $\omega\in \Omega.$

 SPOAFD (identical with SPOAFD2 in the terminology of \cite{QianSAFD}), is for a general Bochner space $ \mathcal{N}=L^2(\Omega,L^2({\D})),$ where the space $L^2({\D})$ has a BVC dictionary. We have the same result except that the TM system $\{B_k\}$ is replaced by the orthonormal system $\{E_k\},$ as composed in (\ref{GSmultiple}), using the multiple kernels $\tilde{K}_{q_k},$  where $q_k$ are selected according to the \emph{Stochastic Pre-orthogonal Maximal Selection Principle} (SPOMSP):
 \begin{eqnarray}\label{MSPD2} q_k=\arg \max \{\mathbb{E}_{\omega} |\langle f_{\omega},E_k^q\rangle |^2\ |\ q\in {\D}\}.\end{eqnarray}
  In the case, as proved in \cite{QianSAFD}, there holds for a.s. $\omega\in \Omega$
 \begin{eqnarray}\label{hold} f(x,\omega)\stackrel{\mathcal{N}}{=}\sum_{k=1}^\infty \langle f_\omega,E_k\rangle E_k(x).\end{eqnarray}

 The power of SPOAFD is that the selected optimal parameters $q_k$'s generate an orthonormal system that uniformly suits for a.s. $\omega,$ and gives rise to, as the Karhunen-Lo\`eve (KL) decomposition, optimal convergence rate.  SPOAFD, as an extension of SAFD, can be associated with any BVC dictionary. As for the deterministic case studied in \cite{Qu-Sps-represt-dirac}, such flexibility makes SPOAFD to be a convenient tool to solve Dirichlet boundary value and Cauchy initial value problems with random data (\cite{YangF}): We take the dictionaries of the shifted and dilated Poisson and heat kernels, respectively, and make use the lifting up technology based on the semigroup properties of the kernels, to get the solutions of the random boundary and initial value problems. The same method can be used to solve other types deterministic and random boundary and initial value problems (\cite{LLQQ}). SAFD and SPOAFD are considerable convenient in the computation respect, as they rely only on the covariance function but not on the eigenvalues and eigenfunctions of the integral operator having the covariance function as its kernel. On the other hand, the widely used and important KL decomposition methodology crucially relies on the eigenvalues and eigenfunctions.\\

 \noindent (g) Stochastic $n$-Best AFD: The related \emph{stochastic $n$-best} (SnB) approximation problems are first formulated and studied in \cite{QQD} (specially for the stochastic complex Hardy spaces) and further in \cite{Qiannbest} (for a wide class of stochastic reproducing kernel Hilbert spaces).
 The general formulation in \cite{Qiannbest} is as follows. For any  $n$-tuple ${\bf p}=(p_1,\cdots,p_n) \in
 {\D}^n,$ where multiplicity is allowed, there exists an $n$-orthonormal system $\{E_k^{{\bf p}}\}_{k=1}^n,$ generated by the corresponding possibly multiple kernels
 $\{\tilde{K}_{p_k}\}_{k=1}^n$ through the GS process. The
 associated objective function to be maximized is
 \begin{eqnarray}\label{refer}A(f,{\bf p})=\mathbb{E}_{\omega}\left(\sum_{k=1}^n |\langle f_\omega,E^{\bf p}_k\rangle|^2\right).\end{eqnarray}
 In other words, the SnB problem
 amounts to finding ${\bf q}\in {\D}^n$ such that
 \[A(f,{\bf q})=\sup \{A(f,{\bf p})\ |\ {\bf p}\in {\D}^n\}.\]
 The just formulated SnB is so far the state of the art among the existing variations of the AFD type models.   The main goal of the present paper is to compare SPOAFD and SnB  with the KL decomposition method for decomposing stochastic processes, the latter being methodology-wise of the same kind as AFD in contrast with the Wiener chaos one.\\

 In \S 2 we give an account of the KL decomposition in relation to the concerns of this study. We prove some new results in relation to certain nesting aspects of the naturally involved RKHSs as subspaces of $L^2(\partial {\D})$. In \S 3 we prove some convergence rate results, analyze and compare the AFD type and the KL decompositions, and specify, in both the theory and computational aspects, respective particular properties of the two types.  Illustrative experiments are contained in \S 4. Conclusions are drawn in \S 5.

 \section{An account to the KL decomposition}

\def\T{\mathcal T}
\def\C{\mathcal C}
 Either of an orthonormal basis of $L^2(\Omega)$ or one in $L^2(\partial {\D})$ induces a decomposition of $f\in L^2(\Omega, L^2(\partial {\D}).$ The Fourier-Hermite expansion and the KL decomposition correspond to a basis in $L^2(\Omega)$ and $L^2(\partial {\D})),$ respectively. Since KL is of the same type as the AFD ones, and they treat the same types of problems, in this study we restrict ourselves to only analyze the KL and the AFD types decompositions.\\

The following material of the KL decomposition is standard (\cite{LPS}). We are involved with a compact set ${\T}$ of the Euclidean space in the Bochner space $L^2(\Omega,L^2({\T}))$ setting, where ${\T}$ is in place of $\partial {\D}$ of the proceeding context.
Let $f(t,\omega)$ be given in $L^2(\Omega,L^2({\T}))$ and fixed throughout the rest of the paper, and $\mu (t)={\mathbb E}_\omega f(t,\cdot).$ Denote by ${\C}$ the covariance function:
\[ {\C}(s,t)={\mathbb E}_\omega \left[(f(s,\cdot)-\mu(s))(\overline{f}(t,\cdot)-\overline{\mu}(t)))\right],\]
and by $T$ the kernel-integral operator using ${\C}(s,t)$ as kernel, $T: L^2({\T})\to L^2({\T}),$
\[ TF(s)=\int_{\T}{\C}(s,t)F(t)dt.\]

We denote by $R(T)$ the range of the operator $T,$ and by $R(L)$ the range of the operator $L: L^2(\Omega)\to L^2({\T}),$ defined as, for any $g\in L^2(\Omega),$
\[ L(g)(t)=E_\omega (g(f^t-\mu(t)))=\int_{\Omega} [{f}(t,\omega)-\mu(t)]g(\omega)d{\mathbb P},\]
where $f^t(\omega)=f(t,\omega).$
It is asserted that $T$ is a Hilbert-Schmidt operator in $L^2({\T}),$ and hence a compact operator. The kernel function $\C$ is conjugate-symmetric and non-negative. As a consequence, $T$ has a sequence of positive eigenvalues $\lambda_1\ge\lambda_2\ge\cdots\ge\lambda_n\ge\cdots>0,$ and correspondingly a sequence of orthonormal eigenfunctions $\phi_1,\cdots,\phi_n,\cdots, T\phi_k=\lambda_k\phi_k.$ If ${\rm span}\{ \phi_1,\cdots,\phi_n,\cdots\}\ne L^2({\T}),$ or, equivalently, $R(T)\ne L^2({\T}),$ a supplementary orthonormal system (corresponding to the zero eigenvalue) may be added to form a complete basis system, called a KL basis, still denoted as $\{\phi_k\}$ with the property $T\phi_k=\lambda_k\phi_k,$ where $\lambda_k$ now may be zero. There holds
 \begin{eqnarray}\label{consequence} {\C}(s,t)=\sum_{k=1}^\infty \lambda_k\phi_k(s)\overline{\phi}_k(t).\end{eqnarray}
 When $\C$ is continuous in ${\T}\times {\T},$ all the $\phi_k$ are continuous, and the above convergence is uniform.
 The originally given $f(\cdot,\omega),$ a.s. belonging to $L^2({\T}),$ has the so called KL decomposition for a.s. $\omega\in \Omega:$
\begin{eqnarray}\label{expansion}f(t,\omega)-\mu(t)\stackrel{\mathcal{N}}{=}
\sum_{k=1}^\infty \sqrt{\lambda_k}\phi_k(t)\xi_k(\omega)\stackrel{\mathcal{N}}{=}
\lim_{n\to\infty}S_n(t,\omega),\end{eqnarray}
where $S_n(t,\omega)\triangleq\sum_{k=1}^n\sqrt{\lambda_k}\phi_k(t)\xi_k(\omega), S_0=0,$ and for the non-zero $\lambda_k,$
\[ \xi_k(\omega)=\frac{1}{\sqrt{\lambda_k}}\langle f_\omega-\mu,\phi_k\rangle_{L^2({\T})}.\]
The random variables are uncorrelated, zero mean, and of unit variance. If the process is Gaussian, then $\xi_k\sim N(0,1)$ iid.

Since ${\C}(s,t)$ is conjugate-symmetric and non-negative, by the Moor-Aronszajn Theorem, it uniquely determines a RKHS. The RKHS is formulated as follows. Let $\alpha_k=1$ if $\lambda_k>0;$ and $\alpha_k=0$ if $\lambda_k=0.$  Define
\[ H_{\C}=\{F\in L^2({\T})\ |\ \|F\|^2_{H_{\C}}=\sum_{k=1}^\infty \alpha_k\frac{|\langle F,\phi_k\rangle|^2}{\lambda_k}<\infty\},\]
whose inner product is defined as
\[ \langle F,G\rangle_{H_{\C}} = \sum_{k=1}^\infty \alpha_k\frac{\langle F,\phi_k\rangle
\overline{{\langle G,\phi_k\rangle}}}{\lambda_k},\]
where the role of the $\alpha_k$ makes what $\lambda_k$ appearing in the sum are those being non-zero.

We collect, in the following theorem,
 some fundamental results of the KL decomposition, of which the ones concerning the nesting relations of the Sobolev type RKHSs may be the first time to be noted.

\begin{theorem}\label{RKHS}
~\\
(i) $H_{\C}\subset L^2({\T});$\par
\smallskip
\noindent(ii) ${\rm var}[f(t,\cdot)-S_n(t,\cdot)]={\C}(t,t)-\sum^n_{k=1}\lambda_k\phi^2_k(t),$  $\|f-S_n\|_{\N}=
\sum_{k=n+1}\lambda_k,$\\ $\|{\rm var}f\|_{L^2(\partial {\D})}=\sum_{k=1}^\infty \lambda_k,$ and
$\|S_n\|_{\N}^2=
{\mathbb E}\sum_{k=1}^n |\langle f_\omega-\mu,\phi_k\rangle|^2=\sum_{k=1}^n \langle T\phi_k,\phi_k\rangle=\sum_{k=1}^n\lambda_k;$\par
\smallskip
\noindent(iii)  $H_{\C}$ is the RKHS with reproducing kernel ${\C}(s,t),$ and, in particular,
$\overline{{\rm span}\{{\C}(s,\cdot)\}}=H_{\C},$ where the bar means the closure under the $H_{\C}$ norm;\par
\smallskip
\noindent(iv) The KL basis has the optimality property: For any orthonormal basis $\{\psi_k\}$ of $L^2(\partial {\D})$ and any $n,$ there holds $\sum_{k=1}^n\langle T\psi_k,\psi_k \rangle\leq \sum_{k=1}^n\langle T\phi_k,\phi_k\rangle;$\par
\smallskip
 \noindent(v) If there are only finitely many $\lambda_k$'s being non-zero, then $f_\omega\in H_{\C}$ for a.s. $\omega\in \Omega;$\par
\smallskip
 \noindent(vi) $R(L)=H_{\C}$ in the set-theoretic sense, and moreover, $H_K=H_{\C},$ where $H_K$ is the RKHS over $R(L)$ defined with the ${\mathcal H}$-$H_K$ formulation (as in \cite{Qian2020HHK}). In particular, when taking $L^2(\Omega)$ as $L^2(\partial {\D})$ and $f$ the parameterized Szeg\"o kernel, we have, as isometric spaces, $R(L)=H_K=H_{\C}=H^2({\D}_1),$ latter being the classical Hardy space in the disc;\par
 \smallskip
 \noindent(vii) Under the natural inner product of $R(T)$ as given in  $H_{{\C}_j}$ below, the identity mapping $R(T)\to H_{\C}\to L^2(\partial {\D})$ are bounded imbeddings. $R(T)=H_{\C}=L^2(\partial {\D})$ if and only if there exist only a finite number of non-zero $\lambda_k$'s; and\par
 \smallskip
 \noindent(viii) Set ${\C}_j(s,t)=\sum_{k=1}^\infty \lambda_k^{j+1}\phi_k(s)\overline{\phi}_k(t), j=0, 1,\cdots,$  where ${\C}_0={\C}, $ and
\[ H_{{\C}_j}=\{F\in L^2(\partial {\D})\ :\ \|F\|^2_{H_{{\C}_j}}=\sum_{k=1}^\infty \alpha_k\frac{|\langle F,\phi_k\rangle|^2}{\lambda_k^{j+1}}<\infty\}.\]
Then $H_{{\C}_j}$ are RKHSs with, respectively, the reproducing kernels
${\C}_j,$ and $TH_{{\C}_j}=H_{{\C}_{j+2}}.$
\end{theorem}

The assertions (i) to (iv) are known knowledge that can be found in the basic literature. See, for  instance, \cite{LPS}. For self-containing, as well as for the algorithm concerns, we provide their proofs, in which the proof of (iv) by using the simplex algorithm may be new.  We give detailed proofs for the newly stated results (v) to (viii). \\

\noindent{\bf Proof}. (i) is the validity of the definition of $H_{\C}.$ \\

(ii) follows from uncorrelation of the $\xi_k$'s. See, for instance \cite{LPS}. We now deduce some relations for later use.
 Resulted from the orthonomality of $\phi_k,$ we have $\|S_n\|_{\N}^2={\mathbb E}\sum_{k=1}^n |\langle f_\omega,\phi_k\rangle|^2.$ Then
\begin{eqnarray}\label{deduction}
{\mathbb E}\sum_{k=1}^n |\langle f_\omega-\mu,\phi_k\rangle|^2&=&
{\mathbb E}\sum_{k=1}^n \langle f_\omega-\mu,\phi_k\rangle{\overline{\langle f_\omega-\mu,\phi_k\rangle}}\nonumber \\
&=& \sum_{k=1}^n \int_{\partial {\D}}\int_{\partial {\D}}{\mathbb E}([f(s,\cdot)-\mu]{\overline{[f(t,\cdot)-\mu]}})
{\phi_k}(t){\overline{\phi}_k}(s)dtds\nonumber \\
&=& \sum_{k=1}^n \langle T\phi_k,\phi_k\rangle\nonumber \\
&=&\sum_{k=1}^n\lambda_k.
\end{eqnarray}\\

Now show (iii). By definition, $H_{\C}$ is a Hilbert space. For any fixed $s\in {\T},$  as a consequence of (\ref{consequence}),
 \[\langle{\C}_s,\phi_k\rangle=\lambda_k\phi_k(s).\]
 Hence,
 \[\sum_{k=1}^\infty \alpha_k\frac{|\langle {\C}_s,\phi_k\rangle|^2}{\lambda_k}=\sum_{k=1}^\infty\lambda_k\phi^2_k(s).\]
 Since the $L^1(\partial {\D})$-norm of the last function in $s$ is equal to
 $\sum_{k=1}^\infty\lambda_k<\infty,$ we have $\|{\C}_s\|_{H_{{\C}}}<\infty,$ a.e.
 This implies that for a.e. $s\in \partial {\D},$
 ${\C}_s$ belongs to $H_{\C}.$ Let $F\in H_{\C}$ with $F=\sum_{k=1}^\infty \alpha_k\langle F,\phi_k\rangle \phi_k.$ Then the reproducing property follows:
\begin{eqnarray*} \langle F,{\C}_s\rangle_{H_{\C}} &=&
\sum_{k=1}^\infty \alpha_k\frac{\langle F,\phi_k\rangle \overline{\langle {\C}_s,\phi_k \rangle}}{\lambda_k}\\
&=&\sum_{k=1}^\infty \alpha_k\frac{\langle F,\phi_k\rangle \lambda_k\phi_k(s)}{\lambda_k}\\
&=&F(s).\end{eqnarray*}\\

Now show (iv).  Let $n$ be fixed, and $m\ge n,$
and $T_m$ be the integral operator defined with the kernel
${\C}_m(s,t)=\sum_{j=1}^m \lambda_j\phi_j(s)\overline{\phi}_j(t),$ where $\phi_1,\cdots,\phi_m$ are the first $m$ elements of the entire KL basis $\{\phi_j\}_{j=1}^\infty.$ For any $n$-orthonormal system  $\{\psi_1,\cdots,\psi_n\},$ there holds
\begin{eqnarray*}
\sum_{k=1}^n\langle T_m\psi_k,\psi_k\rangle &=& \sum_{k=1}^n \int_{\T}\int_{\T} {\C}_m(s,t){\psi}_k(t)\overline{\psi}_k(s)dtds\\
&=&\sum_{k=1}^n \int_{\T}\int_{\T}
\sum_{j=1}^m\lambda_j{\phi}_j(s)\overline\phi_j(t){\psi}_k(t)\overline{\psi}_k(s)dtds\\
&=&\sum_{k=1}^n\sum_{j=1}^m\lambda_j|\langle \psi_k, \phi_j\rangle|^2.
\end{eqnarray*}
Denote $c_{kj}=|\langle \psi_k, \phi_j\rangle|^2.$ Now we try to solve the global maximization problem for the linear objective function
\[A=\sum_{k=1}^n\sum_{j=1}^m\lambda_jc_{kj}\]
under the constrain conditions
\begin{equation*}\begin{cases}
\sum_{k=1}^nc_{kj}+\alpha_j=1,\\
 \sum_{j=1}^mc_{kj}+\beta_k=1, \\
c_{kj}\ge 0, \alpha_j\ge 0, \beta_k\ge 0,\\
\lambda_1\ge\cdots\ge\lambda_n\ge\cdots\lambda_m,\\
1\leq k\leq n, 1\leq j\leq m,
\end{cases}
\end{equation*}
where the first two constrain conditions are due to the Bessel inequality.
This is typically a simplex algorithm problem that assumes the greatest possible value at the tip points of the defined region. By the simplex algorithm method, the optimal solution is attainable at $c_{kk}=1, 1\leq k\leq n; c_{kj}=0, k\ne j; \alpha_j=1, n<j\leq m; \beta_k=0, 1\leq k\leq n.$ This solution simply means that $\psi_k=\phi_k, 1\leq k\leq n.$
As a conclusion of the above simplex algorithm solution, for a general $n$-orthonormal system $\{\psi_1,\cdots,\psi_n\},$  there hold
\[ \sum_{k=1}^n\langle T_m\psi_k,\psi_k\rangle\leq \sum_{k=1}^n\lambda_k=\sum_{k=1}^n\langle T_m\phi_k,\phi_k\rangle.\]
Letting $m\to \infty,$ since $\lim_{m\to \infty}T_m\psi_k=T\psi_k, 1\leq k\leq n,$ we obtain the inequality claimed in (iv).\\

 Next, we show (v). For a.s. $\omega\in \Omega$ we have the series expansion (\ref{expansion}).
 Hence,
\begin{eqnarray*}
\|f_\omega-\mu\|^2_{H_{\C}}&=&\sum_{k=1}^\infty \alpha_k\frac{|\langle f_\omega,\phi_k\rangle|^2}{\lambda_k}\\
&=&\sum_{k=1}^\infty \alpha_k|\xi_k|^2.
\end{eqnarray*}
If there are finitely many $\lambda_k$'s being non-zero, then
\[{\mathbb E}_\omega (\sum_{k=1}^\infty \alpha_k|\xi_k|^2)=\sum_{k=1}^\infty \alpha_k<\infty.\]
This implies that in the case for a.s. $\omega \in \Omega$ there holds $\|f_\omega-\mu\|_{H_{\C}}<\infty,$ and thus a.s. $f_\omega-\mu\in H_{\C}.$  \\

 We now show (vi). Let $g(\omega)\in L^2(\Omega).$ Its image under $L$ is $G=Lg,$ that is, by using the KL expansion of $f,$
\[ G(t)=\int_{\Omega}g(\omega)[{f(t,\omega)-\mu(t)}]d{\mathbb P}=\sum_{k=1}^\infty \sqrt{\lambda_k}
\langle g,\overline{\xi}_k\rangle_{L^2(\Omega)} {\phi}_k(t).\]
Now we examine the $H_{\C}$ norm of $G.$ Noting that $\{\overline{\xi}_k\}$ is an orthonormal system in $L^2(\Omega)$, by invoking the Bessel inequality, we have
\[ \|G\|_{H_{\C}}=\sum_{k=1}^\infty\alpha_k|\langle g,\overline{\xi}_k\rangle_{L^2(\Omega)}|^2\leq \alpha_1\|g\|^2_{L^2(\Omega)}<\infty.\]
Therefore, $L(g)\in H_{\C}.$
By the Riesz-Fisher Theorem and the definition of $H_{\C},$ the mapping
$L: {\mathcal H}\to H_{\C}$ is onto. Hence, in the set theoretic sense $R(L)=H_{\C}.$
With the ${\mathcal H}$-$H_K$ formulation of \cite{Qian2020HHK} (also see \cite{Saitoh_RKHS-2016}), through facilitating the specific inner product over $R(L)$ we obtain the RKHS, $H_K,$ for which ${\C}$ is the reproducing kernel. The inner product used there is defined via the equivalent classes in $L^2(\Omega)$ of the form $L^{-1}(G)$ (as set inverse), $G\in H_{\C}.$  Now $H_{\C}$ is also a RKHS with the same reproducing kernel ${\C}.$ The uniqueness part of the Moor-Aronszajn Theorem asserts that the two RKHSs, namely $H_K$ and $H_{\C},$ have to be the same. The latter stands as a realization of the former in terms of the eigenvalues and eigenfunctions of the integral operator $T.$\\

Next, we show (vii), and first show $R(T)\subset H_{\C}.$ Letting $F=\sum_{k=1}^\infty \langle F,\phi_k\rangle \phi_k\in L^2({\T}),$ then
\begin{eqnarray*}
TF(s)=\int_{\T} {\C}(s,t)F(t)dt=\sum_{k=1}^\infty \alpha_k\lambda_k \langle F,\phi_k\rangle\phi_k(s)=\\
=\sum_{k=1}^\infty \tilde{c}_k\phi_k(s), \quad \tilde{c}_k=\alpha_k\lambda_k \langle F,\phi_k\rangle.\end{eqnarray*}
Those coefficients satisfy the condition
\[ \sum_{k=1}^\infty \frac{|\tilde{c}_k|^2}{\lambda_k^2}= \sum_{k=1}^\infty\alpha_k
|\langle F,\phi_k\rangle|^2<\infty.\]
Since $\lambda_k\to 0,$ the above condition implies
\[ \sum_{k=1}^\infty \frac{|\tilde{c}_k|^2}{\lambda_k}<\infty,\]
and hence
 $TF\in H_{\C}.$
By invoking the Riesz-Fisher Theorem, $R(T)$ has its natural inner product in
\[ H_{\C_1}=\{F\in L^2(\partial {\D})\ :\ \sum_{k=1}^\infty
\alpha_k\frac{|\langle F,\phi_k\rangle|^2}{\lambda_k^2}<\infty\}.\]
Below we may identify $R(T)$ with  $H_{\C_1}.$ Since $\lambda_k\to 0,$ being affiliated with the $H_{\C_1}$-norm the identical mapping from $R(T)$ to $H_{\C}$ is a bounded imbedding: In fact,
\[ \|F\|_{H_{\C}}\leq \max\{\lambda_k:\ k=1,2,\cdots\}\|F\|_{H_{\C_1}}.\]
In the chain of the imbeddings $R(T)=H_{\C_1}\subset H_{\C}=R(L)\subset L^2(\partial {\D})$ any of the two spaces, therefore all the spaces, are identical with $L^2(\partial {\D}),$ if and only if there exist only finitely many non-zero $\lambda_k$'s, for, otherwise, $0\ne \lambda_k\to 0,$ then
$\sum_{k=1}^\infty \sqrt{\lambda_k}\phi_k\in L^2(\partial {\D})\setminus H_{\C},$ and $\sum_{k=1}^\infty \lambda_k\phi_k\in H_{\C}\setminus R(T).$ So (vii) is proved.\\

 We finally prove (viii). Obviously, the norm $\|F\|^2_{H_{{\C}_j}}$ is associated with the inner product \[ \langle F,G \rangle_{H_{{\C}_j}}=\sum_{k=1}^\infty \alpha_k\frac{\langle F,\phi_k\rangle \overline{{\langle G,\phi_k\rangle}} }{\lambda_k^{j+1}}.\]
 Under this inner product the reproducing kernel property is verified
 \[ \langle F,({\C}_j)_s\rangle_{H_{{\C}_j}}=\sum_{k=1}^\infty \alpha_k\frac{\langle F,\phi_k\rangle \overline{\langle ({\C}_j)_s,\phi_k\rangle}} {\lambda_k^{j+1}}=\sum_{k=1}^\infty \alpha_k\frac{\langle F,\phi_k\rangle \lambda_k^{j+1}\phi_k(s)}{\lambda_k^{j+1}}=F(s).\]
 To verify the last statement of (viii), as $TF=\sum_{k=1}^\infty \lambda_k\langle F,\phi_k\rangle \phi_k,$ there follows
 \[\sum_{k=1}^\infty \alpha_k\frac{|\lambda_k\langle F,\phi_k\rangle|^2}{\lambda_k^{j+3}}=\sum_{k=1}^\infty \alpha_k\frac{|\langle F,\phi_k\rangle|^2}{\lambda_k^{j+1}}<\infty.\]
The proof of the theorem is complete.\\

\begin{remark} We note that, in the proof of (iv) of Theorem \ref{RKHS}, not only in the particular case $\psi_k=\phi_k, k=1,\cdots,n,$, this special simplex algorithm problem also assumes its greatest maximal value when ${\rm span}\{\psi_1,\cdots,\psi_n\}=
{\rm span}\{\phi_1,\cdots,\phi_n\}.$ In fact, in such cases, $c_{kj}=|\langle \psi_k,\phi_j\rangle|^2 =0, k\leq n<j\leq m,$ and the Bessel inequality becomes the Placherel identity, $\sum_{k=1}^nc_{kj}=\sum_{j=1}^nc_{kj}=1,$ that implies, continuing with the above deduction,
\begin{eqnarray*}
\sum_{k=1}^n\langle T_m\psi_k,\psi_k\rangle =\sum_{k=1}^n\sum_{j=1}^m\lambda_j|\langle \psi_k, \phi_j\rangle|^2=\sum_{j=1}^n\lambda_j\sum_{k=1}^n|\langle \psi_k, \phi_j\rangle|^2=\sum_{j=1}^n\lambda_j.
\end{eqnarray*}
\end{remark}

\section{The SPOAFD methods in comparison with KL}

\subsection{The algorithms of SAFD, SPOAFD and SnB}
We will give the computational details of SPOAFD. SAFD is just
 SPOAFD restricted to the Hardy $H^2({\D}_1)$ or the $H^2({\bf C}^+)$ space. We note that only in these two latter cases rational orthogonal systems as compositions of Szeg\"o kernels and Blaschke products are available. Let $\{E_q\}$ be a BVC dictionary. We are to find, inductively, $q_1,\cdots,q_k,\cdots,$ such that with the notation of \S 1,
\[ q_k=\arg \max \{{\mathbb E}_\omega |\langle f_\omega-\mu,E^q_k\rangle|^2\ |\ q\in {\D}\}.\]
As in the proof of (ii) of Theorem \ref{RKHS}, the quantity ${\mathbb E}_\omega |\langle f_\omega,E^q_k\rangle|^2$ may be computed, for each $n,$ by
\begin{eqnarray}\label{corresponding}
{\mathbb E}_\omega |\langle f_\omega-\mu,E^q_k\rangle|^2&=&
\int_{\partial {\D}}\int_{\partial {\D}}{\C}(s,t){E^q_k}(t)\overline{E}^q_k(s)dtds.
\end{eqnarray}
With the above expression in terms of the covariance function ${\C}$ one can actually work out, with personal computer, the optimal parameters needed by the SPOAFD expansion without information of the eigenvalues and eigenfunctions of the covariance kernel integral operator.  In contrast, the eigenvalue and the eigenfunction infirmation are crucially required in order to execute the KL expansion.  However, for a general kernel operator the exact information of its eigenvalues and eigenfunctions are by no means easy, and  nor practical. Practically one can only get, by using linear algebra based on samplings, numerical approximations of the eigenvalues and eigenfunctions. With SPOAFD, under a sequence of optimally selected parameters according to (\ref{corresponding}) the relation (\ref{hold}) holds. With SnB for a given $n$ the objective function (\ref{corresponding}) is replaced by (\ref{nB}) seeking for an $n$-tuple of parameters ${\bf q}=(q_1,\cdots,q_n)$ that maximizes
\begin{eqnarray}\label{nB}
\sum_{k=1}^n
\int_{\partial {\D}}\int_{\partial {\D}}{\C}(s,t){E^{\bf p}_k}(t){\overline{E}^{\bf p}_k}(s)dtds
\end{eqnarray}
over all ${\bf p}=(p_1,\cdots,p_n)\in {\D}^n.$ To actually find a practical solution to
(\ref{nB}) one can use the cyclic algorithm given by \cite{Qiancyclic}, or its improvement. Based on the main estimation proved in \cite{Qiannbest} and the algorithm for global maxima of Lipschitz continuous functions on compact sets (For instance, as one given in \cite{QDCZ}), one can get a theoretical algorithm preventing snaking into the local maxima. \\

There is a particular type of random fields whose SPOAFD computation may be simpler. If $f(t,\omega)$ has the form $F(t,X)$ where $X$ is a random variable having a probability density function $p(u), u\in U\subset (-\infty,\infty),$   then (\ref{corresponding}) may be computed as
\[ {\mathbb E}_\omega |\langle f_\omega-\mu,E^q_n\rangle|^2 =
\int_U \int_{\partial {\D}} |F(t,u)E^q_n(t)|^2p(u) dt du.\]
See \cite{YangF} for concrete examples. In such situation for the objective function $(\ref{nB})$ for SnB there exists a similar formula.\\

\subsection{Optimality of KL over SPOAFD and SnB}

The assertion (iv) of Theorem \ref{RKHS} is valid to any orthonomal basis or system $\{\psi_k\},$ and, in particular, valid to the complex-valued orthonormal system $\{E_k\}$ obtained from optimally selected  ${\bf q}=(q_1,\cdots,q_n,\cdots)$ under the SPOAFD or SnB maximal selection principle. We, in fact, have,
\begin{eqnarray}\label{dedu}
{\mathbb E}\sum_{k=1}^n |\langle f_\omega-\mu,E_k\rangle|^2
&=& \sum_{k=1}^n \langle TE_k,E_k\rangle\nonumber \\
&\leq& \sum_{k=1}^n \langle T\phi_k,\phi_k\rangle\nonumber \\
&=&\sum_{k=1}^n \lambda_k.
\end{eqnarray}
We specially note that the functions $E_k$'s are in general not eigenfunctions of the operator $T.$
In spite of the optimality proved in (iv), experiments in \S 4 show, as well as the analysis in the proof of (iv) of Theorem \ref{RKHS} exhibits, when $n$ is large the efficiency of SnB is almost the same as KL. On the other hand, the efficiencies of SPOAFD and SAFD, are very close to SnB. The former two can greatly reduce computation complexity and computer consumes.\\

We also take the opportunity to note that in the POAFD algorithm $|\langle f,E_k\rangle|$ are not necessarily in the descending order. Consider $e_1, e_2$ being two dictionary elements and $f$ is in the span of $e_1, e_2.$ Assume that $\|f_{e_1}\|>|\langle f,e_1\rangle|> |\langle f,e_2\rangle|,$ where $f_{e_1}$ is denoted as the projection of $f$ into the subspace perpendicular to $e_1.$ Then,
$|\langle f,E_2^{e_2}\rangle |=|\langle f_{e_1},E_2^{e_2}\rangle |=\|f_{e_1}\|>|\langle f,e_1\rangle|=|\langle f,E_1\rangle |,$ as claimed.
 \\

\subsection{Convergence Rates}

 To deduce a convergence rate for SPOAFD we adopt a slightly different but more practical formulation. The following maximal selection principle is a little weaker than what is posted in (\ref{MSPD2}), called Weak Stochastic Maximal Selection Principle (WSPOMSP), and the corresponding algorithm is phrased as WSPOAFD with the self-explanatory meaning.

 For $\rho\in (0,1]$ and each $k$, WSPOMSP involves determination of
 a $q_k\in {\D}$ such that
 \begin{eqnarray}\label{WMSPD} \mathbb{E}_{\omega} |\langle f_{\omega}-\mu,E_k^{q_k}\rangle |^2 \ge \rho \sup \{\mathbb{E}_{\omega} |\langle f_{\omega}-\mu,E_k^q\rangle |^2\ |\ q\in {\D}\}.\end{eqnarray}
 If, in particular, $\rho=1,$ then the WSPOMSP reduces to SPOMSP in (\ref{MSPD2}). In some literature there exist more general setting allowing different $\rho_k\in (0,1]$ for different $k,$ we, however, keep it simple and adopt a $\rho\in (0,1]$ uniform in $k.$
Denote
\begin{eqnarray}\label{M} M(\omega)=\inf\{\sum_{l=1}^\infty |c_l(\omega)|\ |\ f_\omega-\mu=\sum_{l=1}^\infty c_l(\omega)E_{q_l}, \forall l, \ q_l\in {\D}\}\end{eqnarray}
and
\[ M_0=\left(\int_\Omega |M(\omega)|^2d\mathbb{P}\right)^{\frac{1}{2}}.\]
The following is an immediate consequence of the convergence rate of the AFD type algorithms. See, for instance, \cite{Q2D}.

\begin{theorem}\label{convergence rate}
Denoting $g_n$ the $n$-standard remainder of the WSPOAFD algorithm,
\[ g_n(x,\omega)=f(x,\omega)-\mu (x)-\sum_{k=1}^{n-1}
\langle {(f_\omega-\mu (x))}_k,E_{q_k}\rangle E_k(x),\]
there holds the estimation
\[\|g_n\|_{\N}\leq \frac{M_0}{\rho \sqrt{n}}.\]
\end{theorem}
Below we include the one-line proof.\\

\noindent{\bf Proof}.
By invoking the deterministic case result, Theorem 3.3 of \cite{Q2D}, there holds, for a.s. $\omega\in \Omega,$
\[ \|{(g_\omega)}_n\|_{L^2({\D})}\leq \frac{M(\omega)}{\rho \sqrt{n}}.\]
By taking the square-norm of the probability space to both sides we obtain the claimed estimation.\\

For the deterministic case there holds $M(\omega)\equiv M_0,$ and the estimation reduces to the WPOAFD one. For $\rho=1$ the estimation reduces to the one for SPOAFD. \\

Since the one by one parameters selection model in SPOAFD is surely less optimal than SnB, we obtain that the convergence rate for SnB is at least the same as that for SPOAFD given by Theorem
\ref{convergence rate}. Note that if the used dictionary is itself an orthonormal system, then the above used convergence rate in the deterministic case quoted from \cite{Q2D} coincides with that for $p=1$ in \cite{DT}.

 \begin{theorem}\label{KL convergence rate}
 Let $n$ be a fixed positive integer. For a given $f\in L^2(\Omega, L^2(\partial {\D})),$ its KL $n$-partial sum expansion, that is the projection of $f$ into the span by the first $n$ eigenfunctions of the integral operator with the kernel $\C,$  coincides with its SnB expansion, where the latter uses the dictionary consisting of the KL basis $\{\phi_l\}_{l=1}^\infty.$ In the case
 \begin{eqnarray}\label{RHS}
 \|f-\mu-\sum_{l=1}^n\sqrt{\lambda_k}\phi_l\xi_l\|_{\N} \leq
\frac{\|\sum_{l=1}^\infty |\langle f_\omega-\mu,\phi_l\rangle|\|_{L^2(\Omega)}}{\sqrt{n}}.\end{eqnarray}
 \end{theorem}
\noindent {\bf Proof}.
 Due to uniqueness of expansion of $f_\omega$ in the basis $\{\phi_k\},$ the infimum in (\ref{M}) reduces to
\begin{eqnarray}\label{uuse}
 M(\omega)=\sum_{l=1}^\infty | \langle f_\omega-\mu,\phi_l\rangle |.
 \end{eqnarray}
 As a consequence of the optimality property of the KL basis, the $n$-partial sum expanded by the first $n$ eigenfunctions of $f\in L^2(\Omega, L^2(\partial {\D}))$ is identical with its SnB for the dictionary $\{\phi_l\}_{l=1}^\infty.$ There holds, owing to (\ref{uuse}) and Theorem \ref{convergence rate},
\[ \|f-\mu-\sum_{k=1}^n\sqrt{\lambda_k}\phi_k\xi_k\|_{\N} \leq
\frac{\|\sum_{k=1}^\infty |\langle f_\omega-\mu,\phi_k\rangle|\|_{L^2(\Omega)}}{\sqrt{n}},\]
The proof is complete.\\

We note that the convergence rate of the greedy algorithm type expansions as given in Theorem \ref{convergence rate} is very coarse. In the KL expansion case in general we
have the actual quantity of the energy of the tail being
\[ \left(\sum_{k=n+1}^\infty\lambda_k\right)^{\frac{1}{2}}<\infty.\]
In contrast, the right-hand-side of (\ref{RHS}) is
\[  \frac{\|\sum_{k=n+1}^\infty\sqrt{\lambda_k}|\xi_k|\|_{L^2(\Omega)}}{\sqrt{n}},\]
provided that the last quantity is a finite number. \\

On the other hand, for the particular Brownian bridge case, for instance, the convergence rate can be precisely estimated as
\[ {\mathbb E}\left[\|B-S_n\|^2_{L^2[0,1]}\right]=\sum_{j=n+1}^\infty \frac{1}{\pi^2j^2}\sim \frac{1}{\pi^2n}\]
(page 206 of \cite{LPS}), showing that the rate is not in general better than
$O(\frac{1}{\sqrt{n}}).$  \\

The nesting structure of the RKHSs studied in  (vii) and (viii) of Theorem (\ref{RKHS}) offers a sequence of finer and finer RKHSs.
One may define the degree $d(F)$ of a function $F\in L^2(\partial {\D})$ as
\[ d(F)=\max\{j\ |\ F\in H_{{\C}_{j}}\setminus H_{{\C}_{j+1}}\}.\]
  A particular $F\in L^2(\partial {\D})$ belongs to all $H_{{\C}_{l}}, l\leq d(F).$ In accordance with the convergence rate results obtained in Theorem \ref{convergence rate} and Theorem \ref{KL convergence rate}, one has discretion to perform an AFD type approximation in a selected space $H_{{\C}_{l}}, l\leq d(F),$ to get a more desirable result.

\subsection{Flexibility of Dictionary Taken by SPOAFD}
For a given random signal $f\in L^2(\Omega,L^2(\partial {\D}))$
 both its KL and its SPOAFD decompositions are
 adaptive to the particular properties of $f.$
 SPOAFD, however, possesses extra adaptivity because the dictionary
 in use can be selected at one's discretion, and especially according to the task taken.
As an example, in \cite{YangF}, we numerically solve the Dirichelet problem of random data
\begin{equation}\label{Equ-1}
\left\{
\begin{aligned}
&\Delta u(x,\omega)=0, ~ \forall~x \in {\D}\subseteq \mathbb{R}^{n+1},~{\rm a.s.}~\omega \in \Omega,\\
&u(x,\omega)=f(x,\omega),~ \ {\rm for \ a.e.}\ ~x\in{\partial  {\D}},~{\rm a.s.}~\omega \in \Omega,
\end{aligned}
\right.
\end{equation}
where $f \in  L^2(\Omega,  L^2(\partial {\D})),~ {\D}=B_1.$ According to the related Hardy space theory the solution $u$ will belong to $L^2(\Omega,  h^2({\D})).$
The most convenient
 dictionary that we use in the case is the one consisting of
 the parameterized Poisson kernels (for the unit ball) $P_{x}(y^{\prime})$ defined as
\begin{equation}\label{Poiss-Ker-B_1}
P_x(y^{\prime})\triangleq P(x,y^{\prime})\triangleq c_n \frac{1-r^2}{|x-y^{\prime}|^n},~
x=rx^{\prime} \in B_1,~  x^{\prime}, y^{\prime} \in \partial B_1,
\end{equation}
The random data $f(x,\omega)$ is efficiently approximated on the boundary by its SPOAFD series in the $\N$ norm: ${\N}=L^2(\Omega,  L^2(\partial {\D})),$
\begin{equation}\label{SPOAFD-h1}
	{f_\omega}(x^{\prime})-\mu(x^{\prime})\stackrel{\mathcal{N}}{=}\sum_{k=1}^{\infty}\langle f_\omega-\mu,E_k\rangle {E_k}(x^{\prime})=\sum_{k=1}^\infty c_k(\omega)\tilde{P}_{x_k}(x^{\prime}), \quad x^{\prime}\in B_1,
	\end{equation}
where $(E_1,\cdots,E_k,\cdots)$ are, consecutively, the GS orthogonalization of
the multiple kernels $(\tilde{P}_{x_1},\cdots,\tilde{P}_{x_k},\cdots), x_k=r_kx_k^{\prime}, k=1,\cdots,n,\cdots,$ and $c_k(\omega)$'s are the associated coefficients converting the two bases. Then based on the semigroup property of the Poisson kernel we can $\lq\lq$lift up" the series expansion on the sphere to get the solution to the Dirichlet problem with the random data, i.e., in the self-explanatory notation,
\begin{equation}\label{solut-Poiss-Integl}
u_{f_\omega-\mu}(x)=\sum_{k=1}^\infty c_k(\omega)\tilde{P}_{x_k}(x), ~\quad x \in B_1.
\end{equation}
The convergence speed of the $n$-partial sums of (\ref{solut-Poiss-Integl}) is the same as that for (\ref{SPOAFD-h1}), being of the rate $O(\frac{1}{\sqrt{n}})$ as proved in \S 3.  By using other bases such as the KL or the Wiener chaos ones the boundary data can also be efficiently expanded but there will be no convenience such as the lifting up associated with the particular equations.
 Some other similar examples are also given in \cite{YangF} and \cite{LLQQ}. We finally note that, not like the KL decomposition, SAFD, SPOAFD and the related SnB are also available in unbounded domains $\D.$

\section{Experiments}
In this section we approximate the Brownian bridge in $[0,2\pi].$ In the following experiments the graphs of the targeted Brownian bridge are made by using the algorithm on page 195 of \cite{LPS}. The KL expansions are according to the formula (5.44) on page 206 of \cite{LPS}. The AFD type methods are based on the covariance of the Brownian bridge, that is
 $$\C(s,t)=\min(s,t)-\frac{st}{2\pi}.$$

\begin{example}(SPOAFD based on the Poisson kernel dictionary)\label{Brownian1}
The experiment is based on 126 sampling points in $[0,2\pi]$ with the uniform spacing $\Delta t \thickapprox 0.05$. As shown in Figure 1,
SPOAFD by using the Poisson kernel dictionary has almost the same effect as that of the KL expansion. At the 125th partial sums both SPOAFD and KL expansions approximately recover the target function. In the local details SPOAFD seems to have visually better results. The relative errors of the two methods are given in Table 1.
\end{example}
\begin{figure}[H]
	\begin{minipage}[c]{0.23\textwidth}
		\centering
		\includegraphics[height=3.5cm,width=3.2cm]{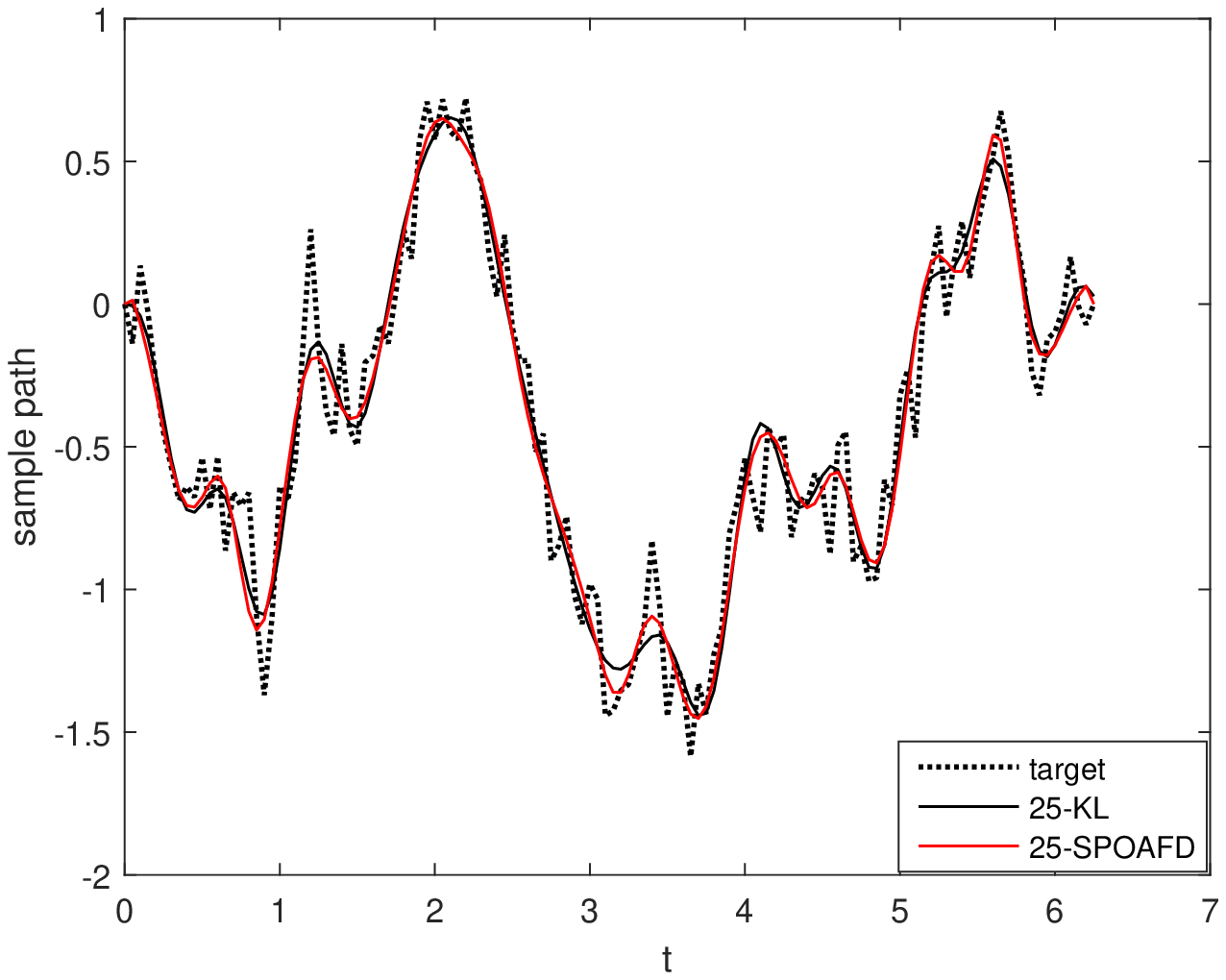}
		\centerline{{\tiny SPOAFD: 25 partial sum}}
	\end{minipage}
	\begin{minipage}[c]{0.23\textwidth}
		\centering
		\includegraphics[height=3.5cm,width=3.2cm]{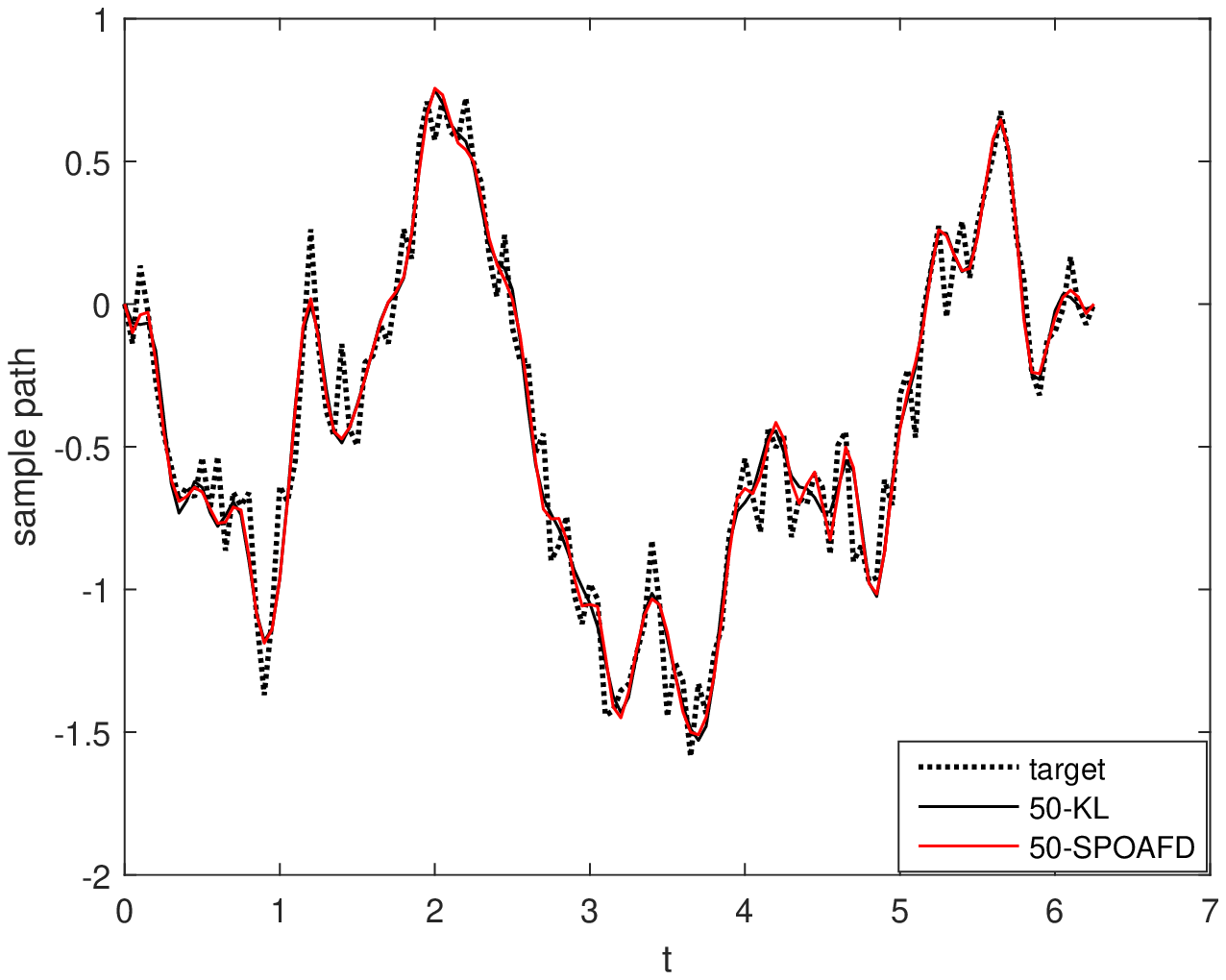}
		\centerline{{\tiny SPOAFD: 50 partial sum}}
	\end{minipage}
	\begin{minipage}[c]{0.23\textwidth}
		\centering
		\includegraphics[height=3.5cm,width=3.2cm]{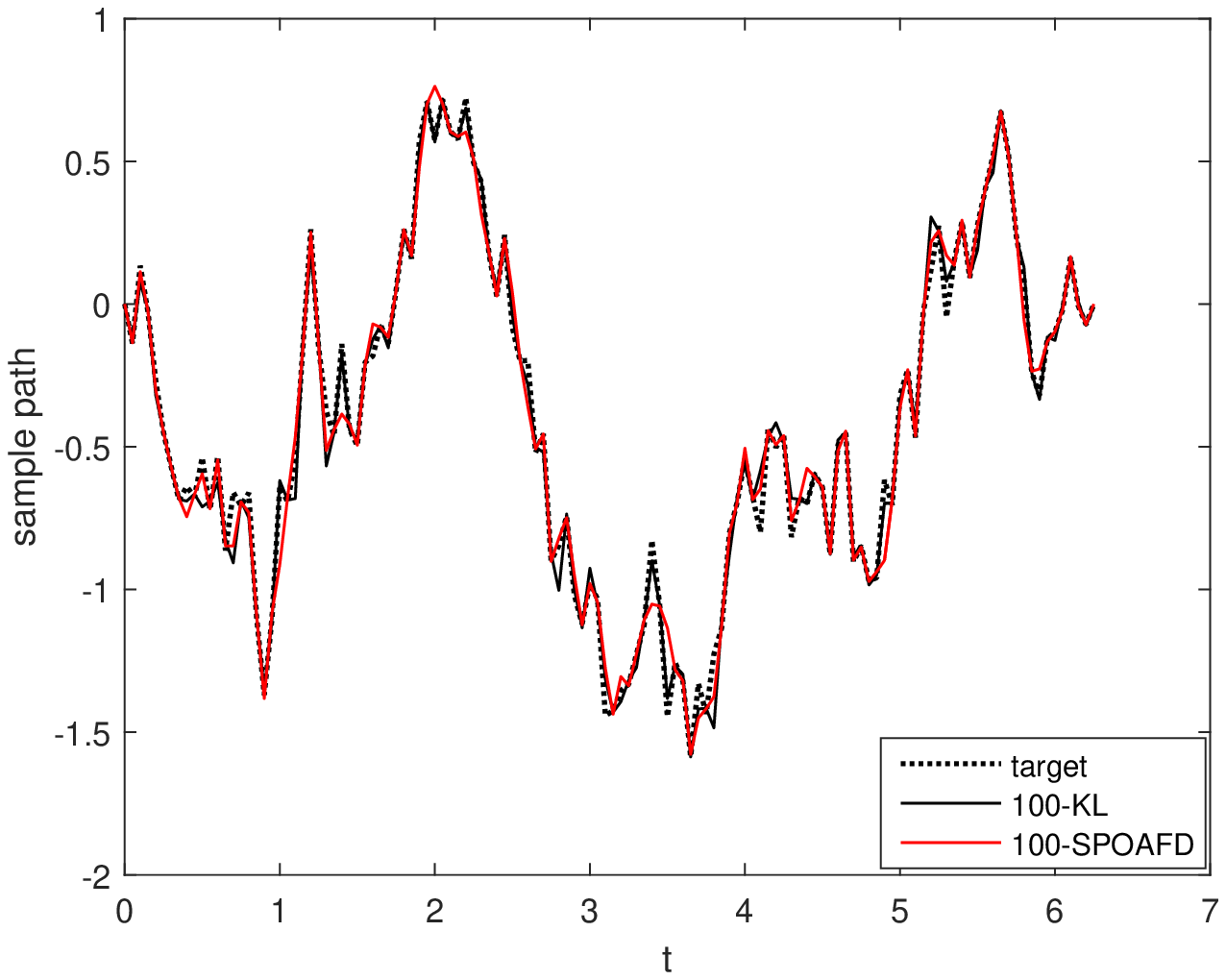}
		\centerline{{\tiny SPOAFD: 100 partial sum}}
	\end{minipage}
	\begin{minipage}[c]{0.23\textwidth}
		\centering
		\includegraphics[height=3.5cm,width=3.2cm]{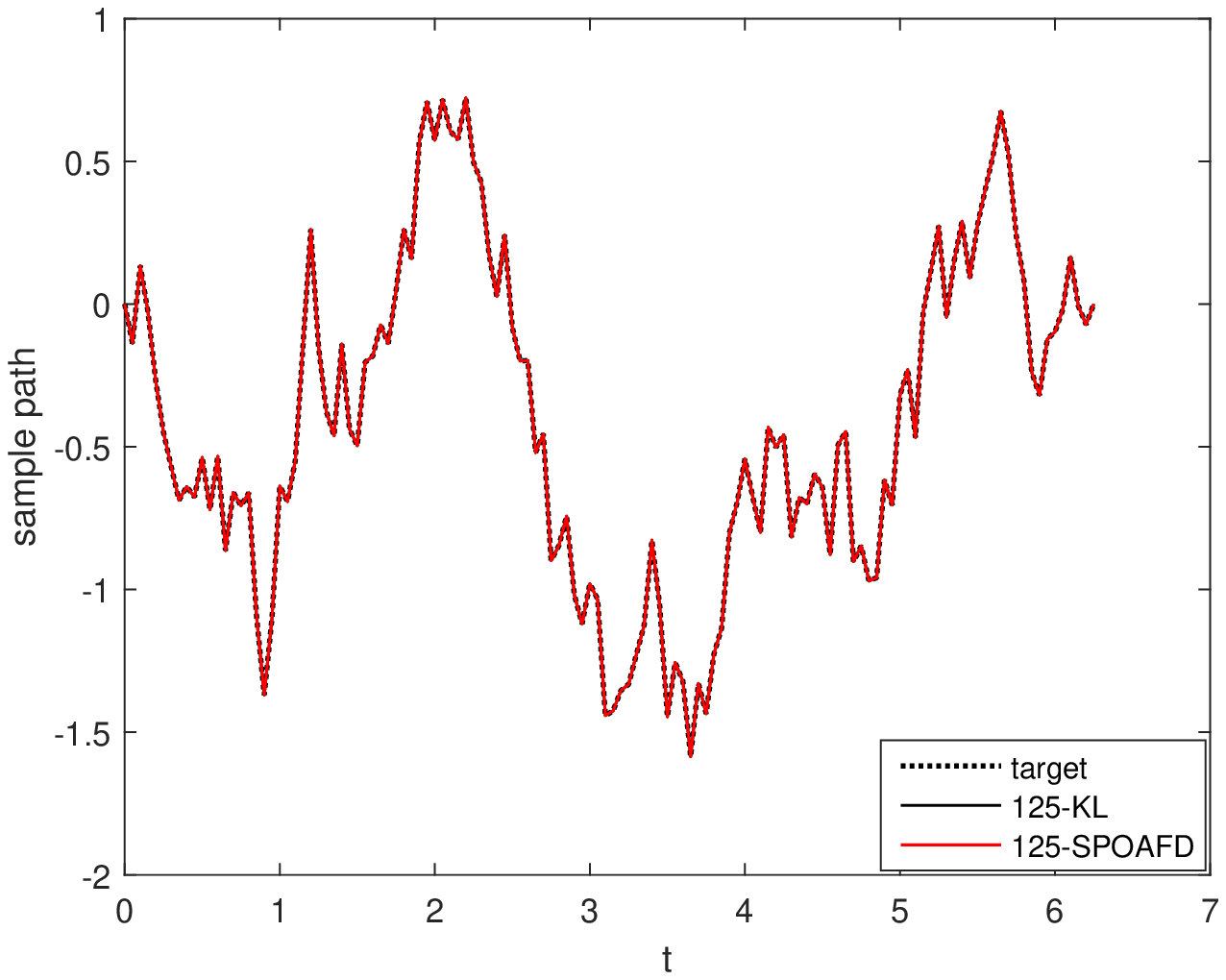}
		\centerline{{\tiny SPOAFD: 125 partial sum}}
	\end{minipage}
	\caption*{Figure 1: sample path I ($126$ points) of Brownian bridge}
\label{figureB1}
\end{figure}
\begin{table}[h]
\scriptsize
  \centering
    \begin{tabular}{c|c|c|c|c}
    \toprule
    \multicolumn{1}{c|}{$n$ partial sum} & \multicolumn{1}{c|}{$n=25$} & \multicolumn{1}{c|}{$n=50$} & \multicolumn{1}{c|}{$n=100$} & \multicolumn{1}{c}{$n=125$} \\
    \midrule
    KL    & 0.0331 & 0.0140 & 0.0021 & $6.1397\times10^{-31}$ \\
    \midrule
    SPOAFD & 0.0298 & 0.0113 & 0.0026 & $1.0984\times10^{-7}$ \\
    \bottomrule
    \end{tabular}%
    \caption*{Table 1: Relative error}
      \label{tableSPOAFD}%
\end{table}%

\begin{example}(SAFD on the Szeg\"o kernel dictionary for the complex Hardy space on the disc space)\label{Brownian2}
We approximate the Brownian bridge by using the KL and the SAFD expansions based on 4096 and 1024 sampling points in $[0,2\pi]$ with the uniform spacing $\Delta t \thickapprox 0.002$ and $\Delta t \thickapprox 0.006,$ respectively. The results are shown in Figure 2a and Figure 2b. SAFD has the convenience of using the TM system that is, in the continuous formulation, orthonormal. Discretely, however, the orthonormal properties are with errors. Hence SAFD requires more sampling points than SPOAFD. The relative errors are given in Table 2a and Table 2b.
\end{example}

\begin{figure}[H]
	\begin{minipage}[c]{0.23\textwidth}
		\centering
		\includegraphics[height=3.5cm,width=3.2cm]{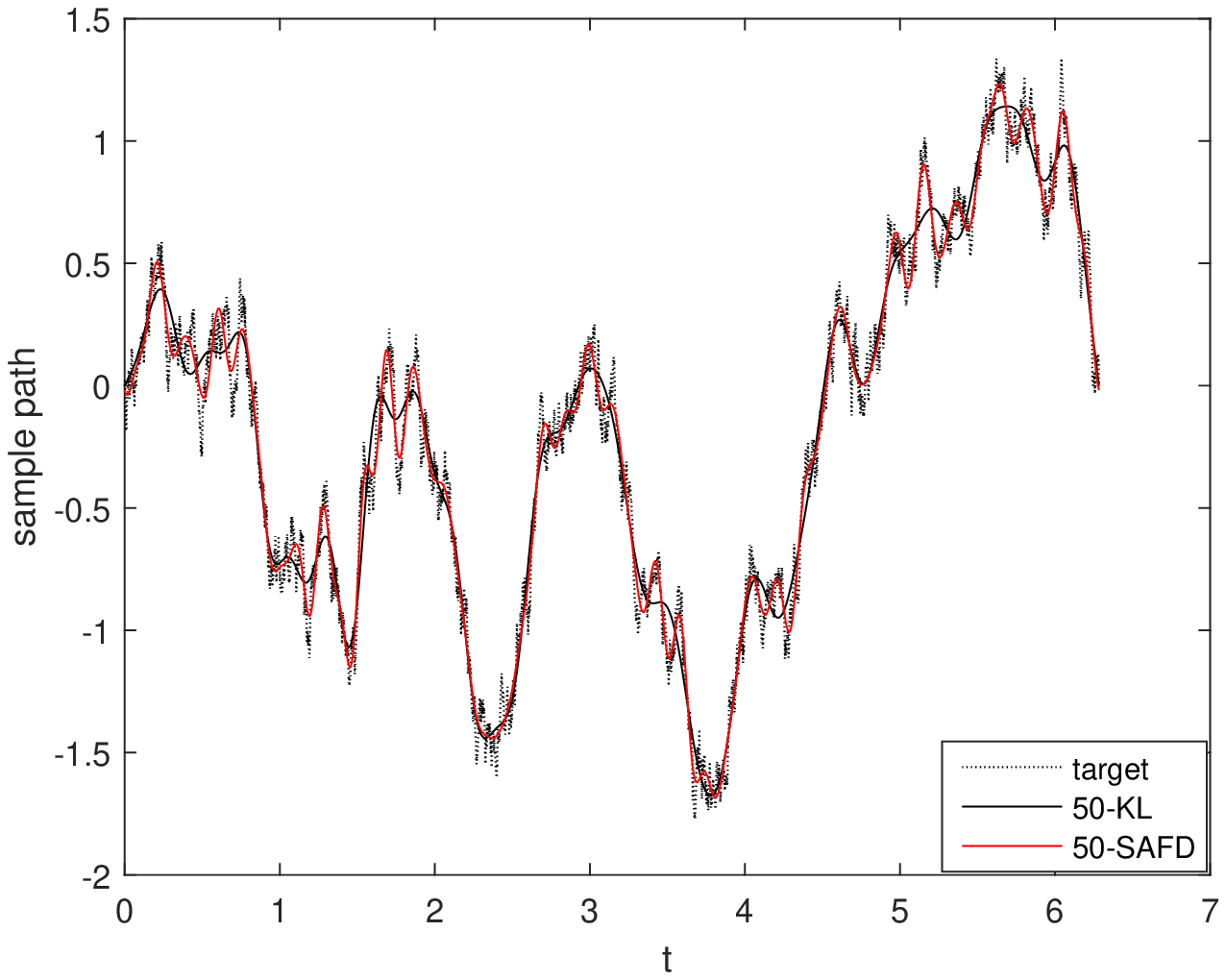}
		\centerline{{\scriptsize SAFD: 50 partial sum}}
	\end{minipage}
	\begin{minipage}[c]{0.23\textwidth}
		\centering
		\includegraphics[height=3.5cm,width=3.2cm]{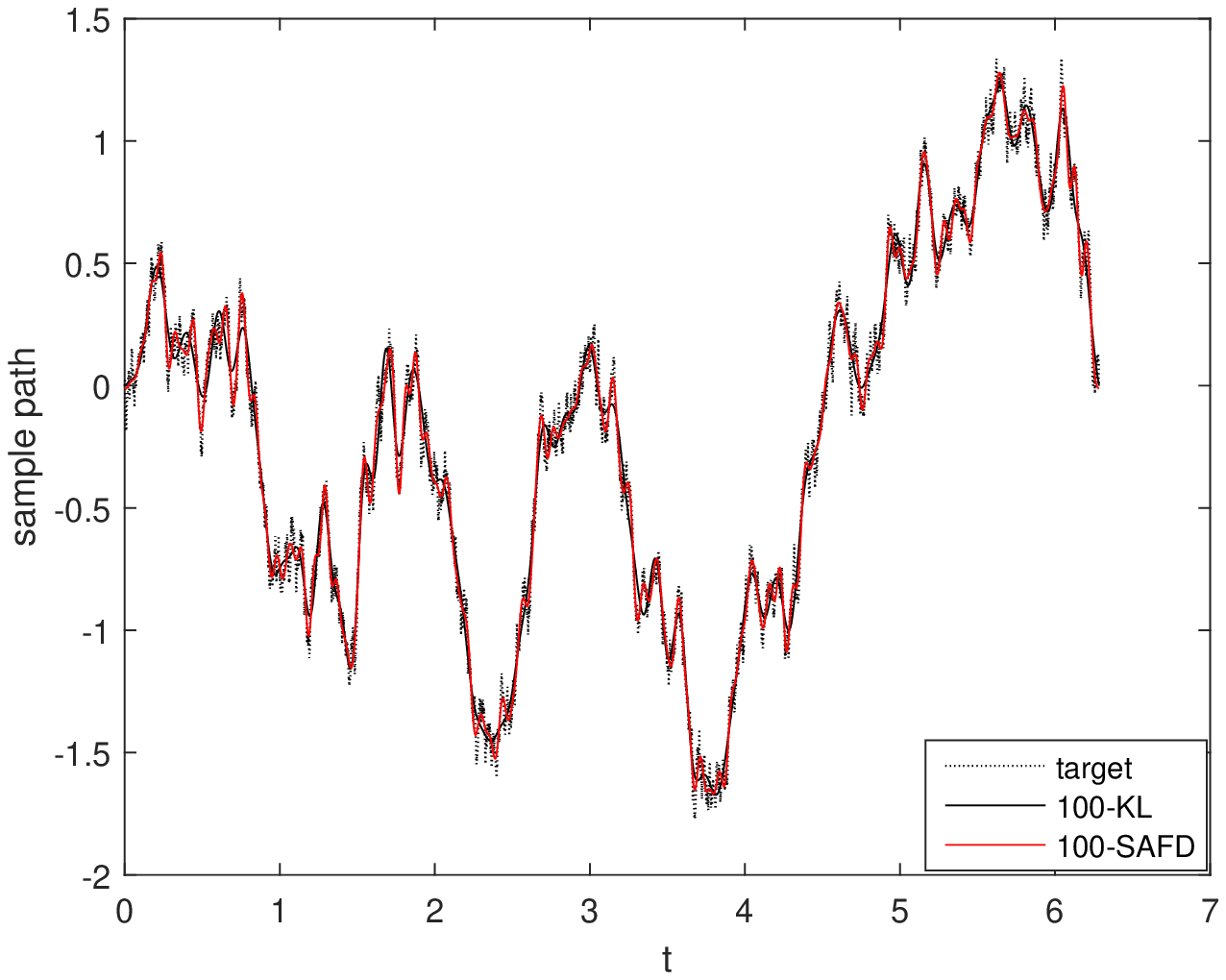}
		\centerline{{\scriptsize SAFD: 100 partial sum}}
	\end{minipage}
	\begin{minipage}[c]{0.23\textwidth}
		\centering
		\includegraphics[height=3.5cm,width=3.2cm]{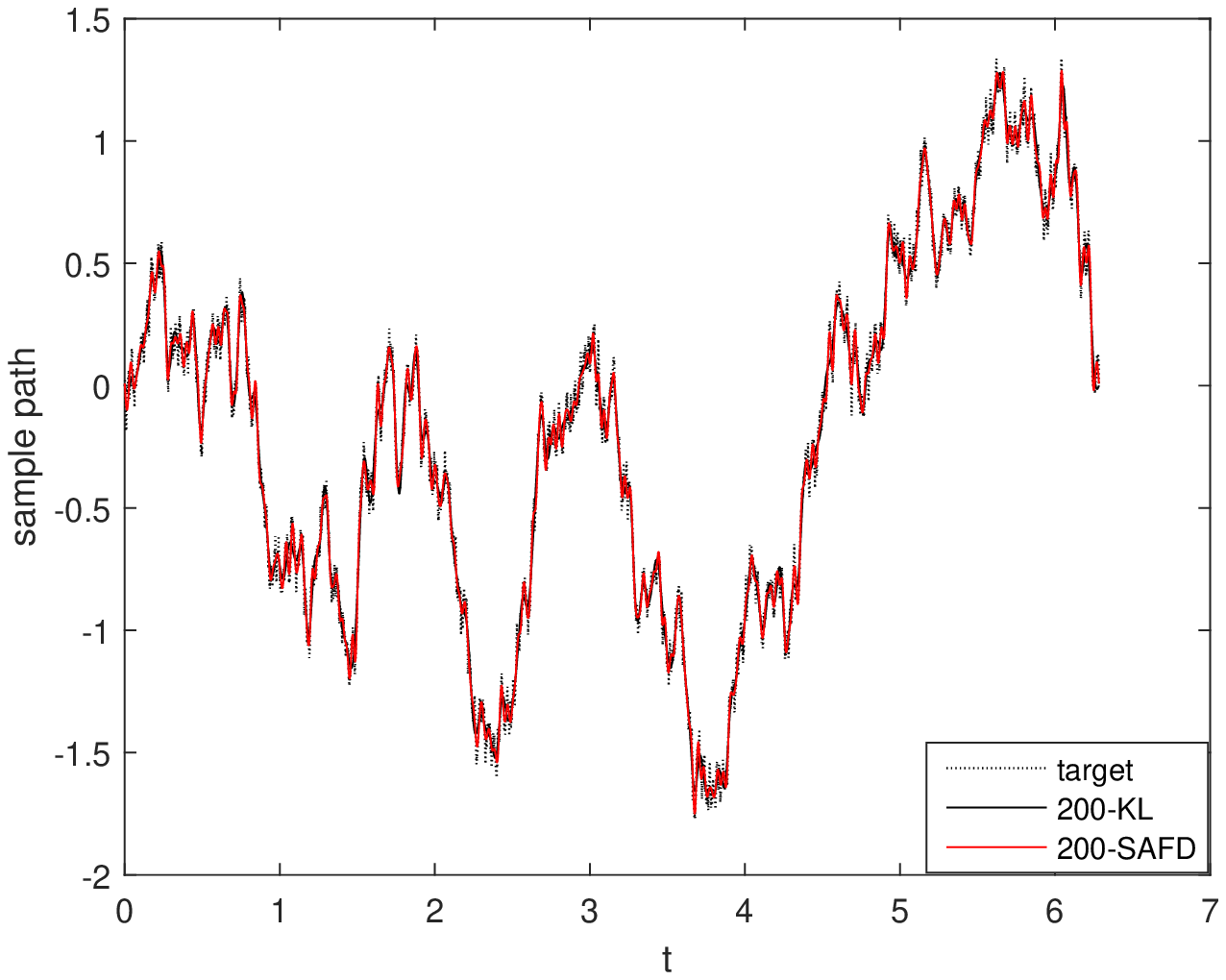}
		\centerline{{\scriptsize SAFD: 200 partial sum}}
	\end{minipage}
	\begin{minipage}[c]{0.23\textwidth}
		\centering
		\includegraphics[height=3.5cm,width=3.2cm]{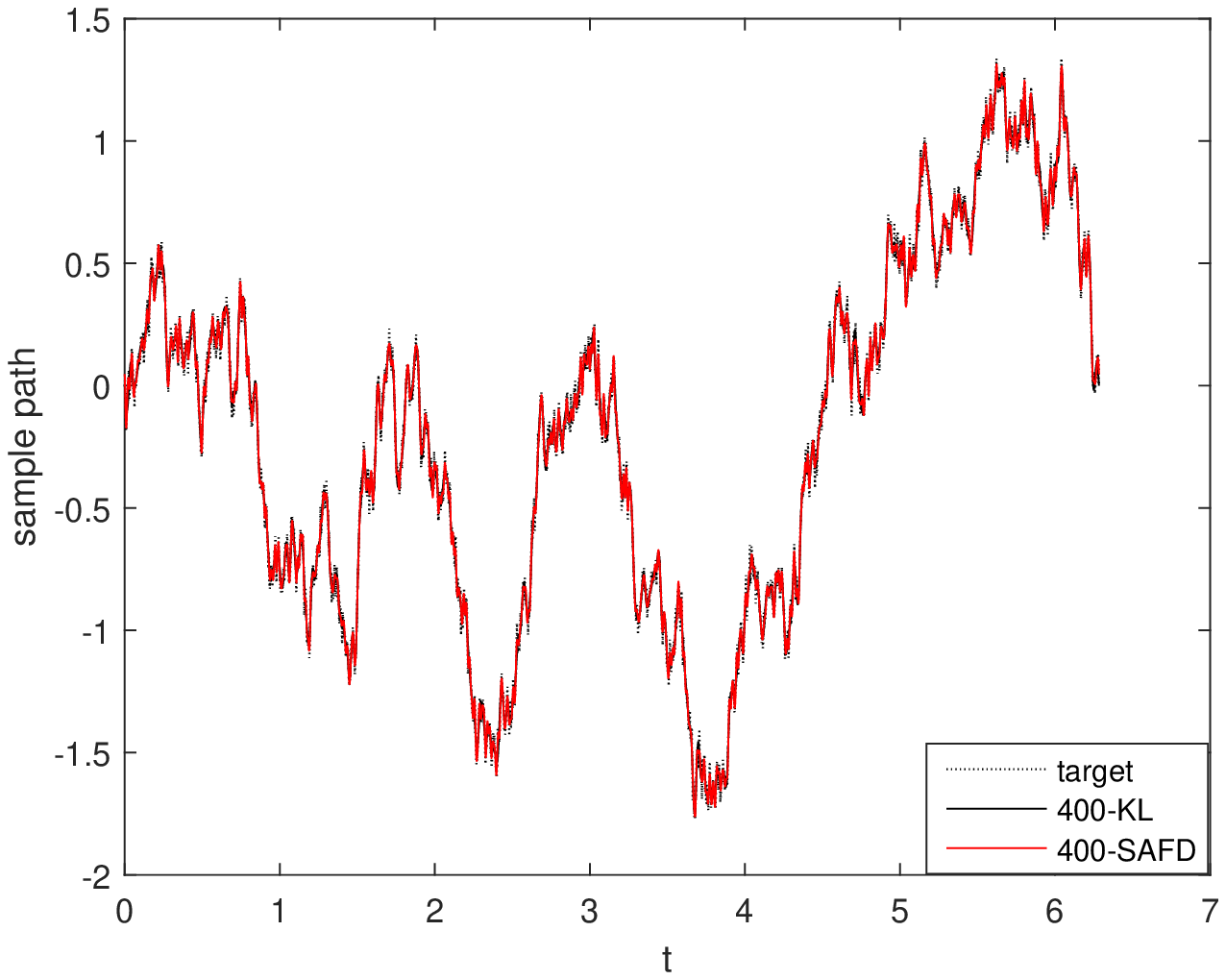}
		\centerline{{\scriptsize SAFD: 400 partial sum}}
	\end{minipage}
	\caption*{Figure 2a: sample path IIa ($2^{12}$ points) of Brownian bridge}
\label{figureB2a}
\end{figure}
\begin{table}[h]
\scriptsize
  \centering
    \begin{tabular}{c|c|c|c|c}
    \toprule
     \multicolumn{1}{c|}{$n$ partial sum} & \multicolumn{1}{c|}{$n=50$} & \multicolumn{1}{c|}{$n=100$} & \multicolumn{1}{c|}{$n=200$} & \multicolumn{1}{c}{$n=400$} \\
    \midrule
    KL    & 0.0237 & 0.0118 & 0.0055 & $0.0026$ \\
    \midrule
    SAFD & 0.0119 & 0.0055 & 0.0026 & $0.0012$ \\
    \bottomrule
    \end{tabular}%
    \caption*{Table 2a: Relative error}
      \label{tableSAFD1}%
\end{table}%

\begin{figure}[H]
	\begin{minipage}[c]{0.23\textwidth}
		\centering
		\includegraphics[height=3.5cm,width=3.2cm]{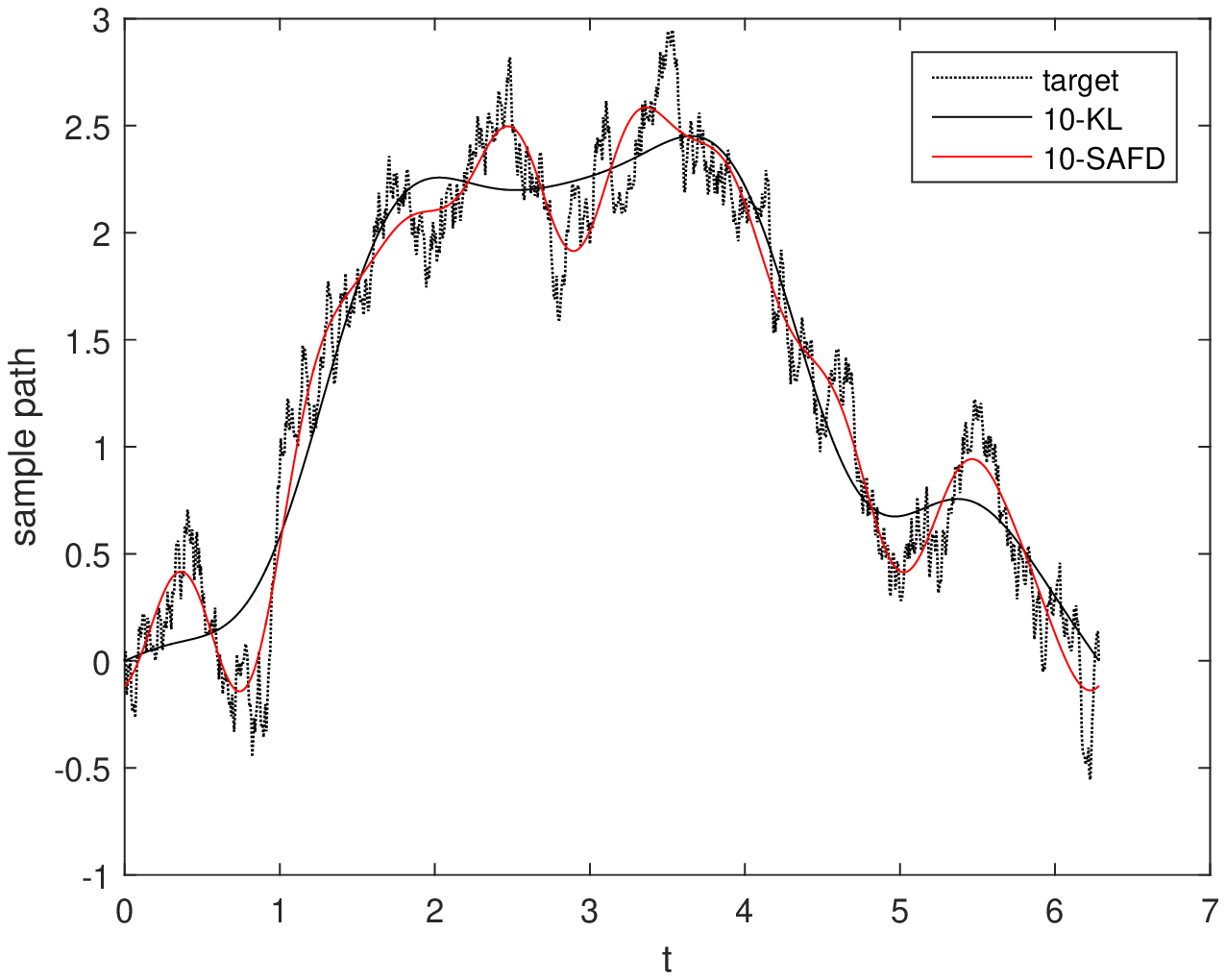}
		\centerline{{\scriptsize SAFD: 10 partial sum}}
	\end{minipage}
	\begin{minipage}[c]{0.23\textwidth}
		\centering
		\includegraphics[height=3.5cm,width=3.2cm]{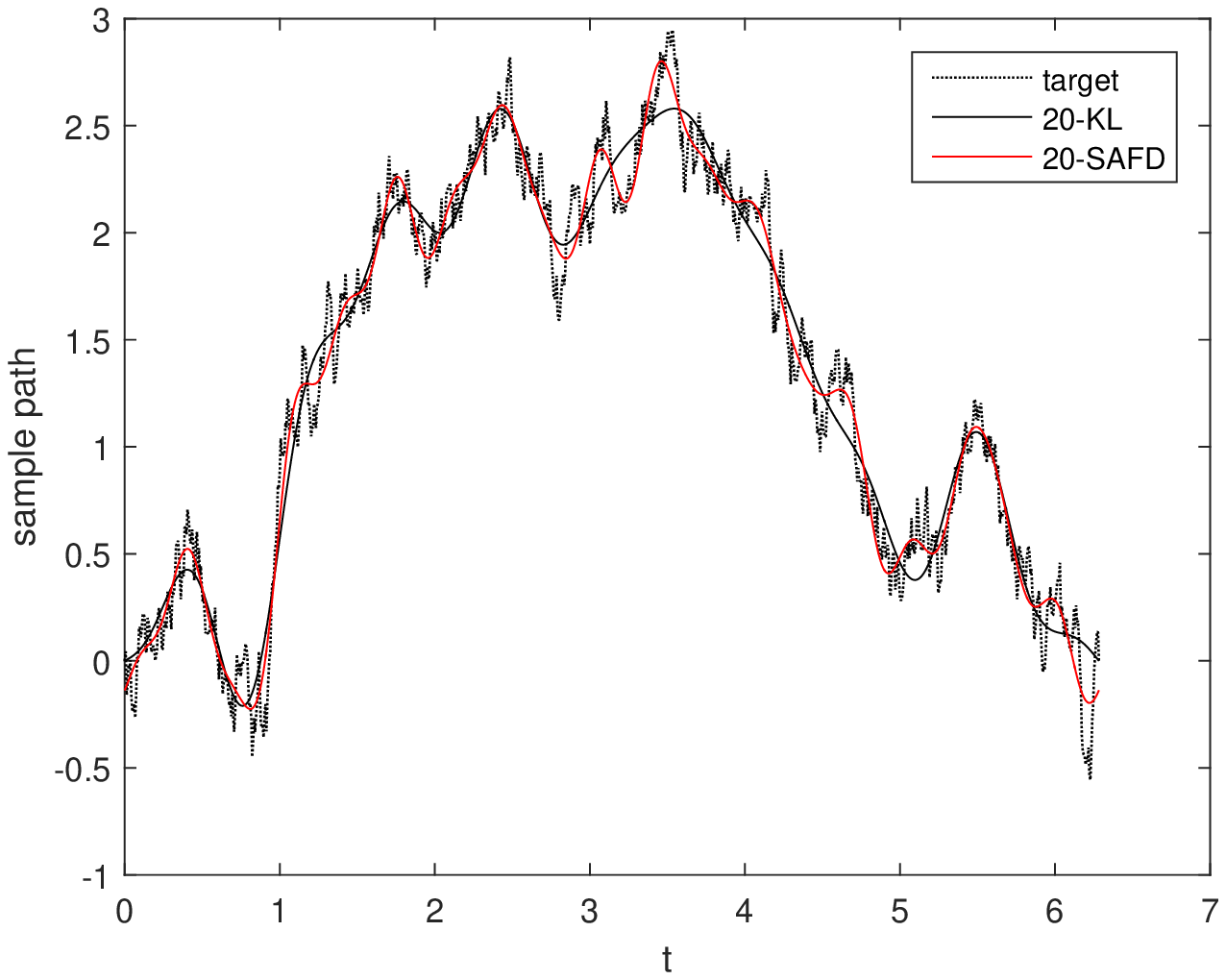}
		\centerline{{\scriptsize SAFD: 20 partial sum}}
	\end{minipage}
	\begin{minipage}[c]{0.23\textwidth}
		\centering
		\includegraphics[height=3.5cm,width=3.2cm]{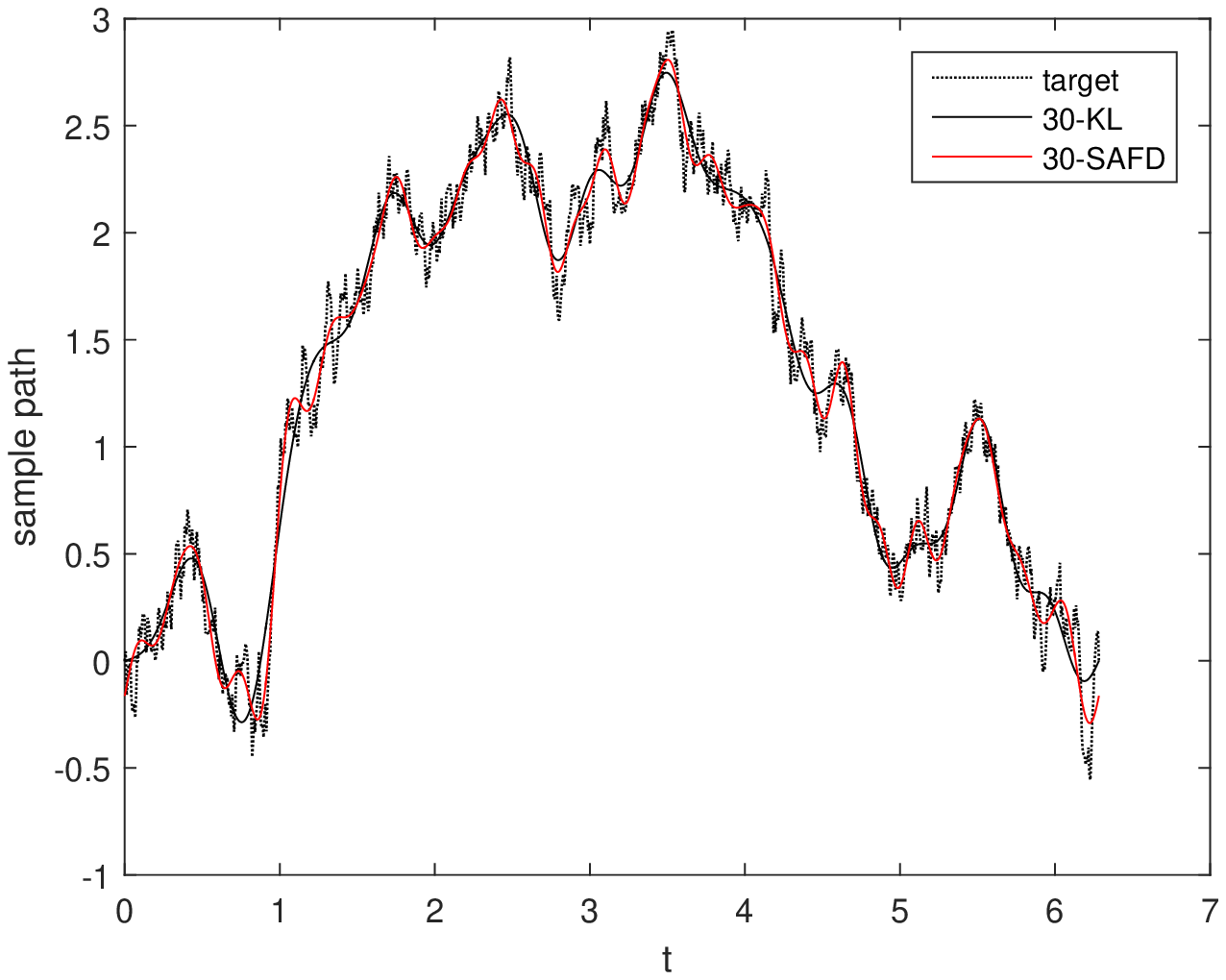}
		\centerline{{\scriptsize SAFD: 30 partial sum}}
	\end{minipage}
	\begin{minipage}[c]{0.23\textwidth}
		\centering
		\includegraphics[height=3.5cm,width=3.2cm]{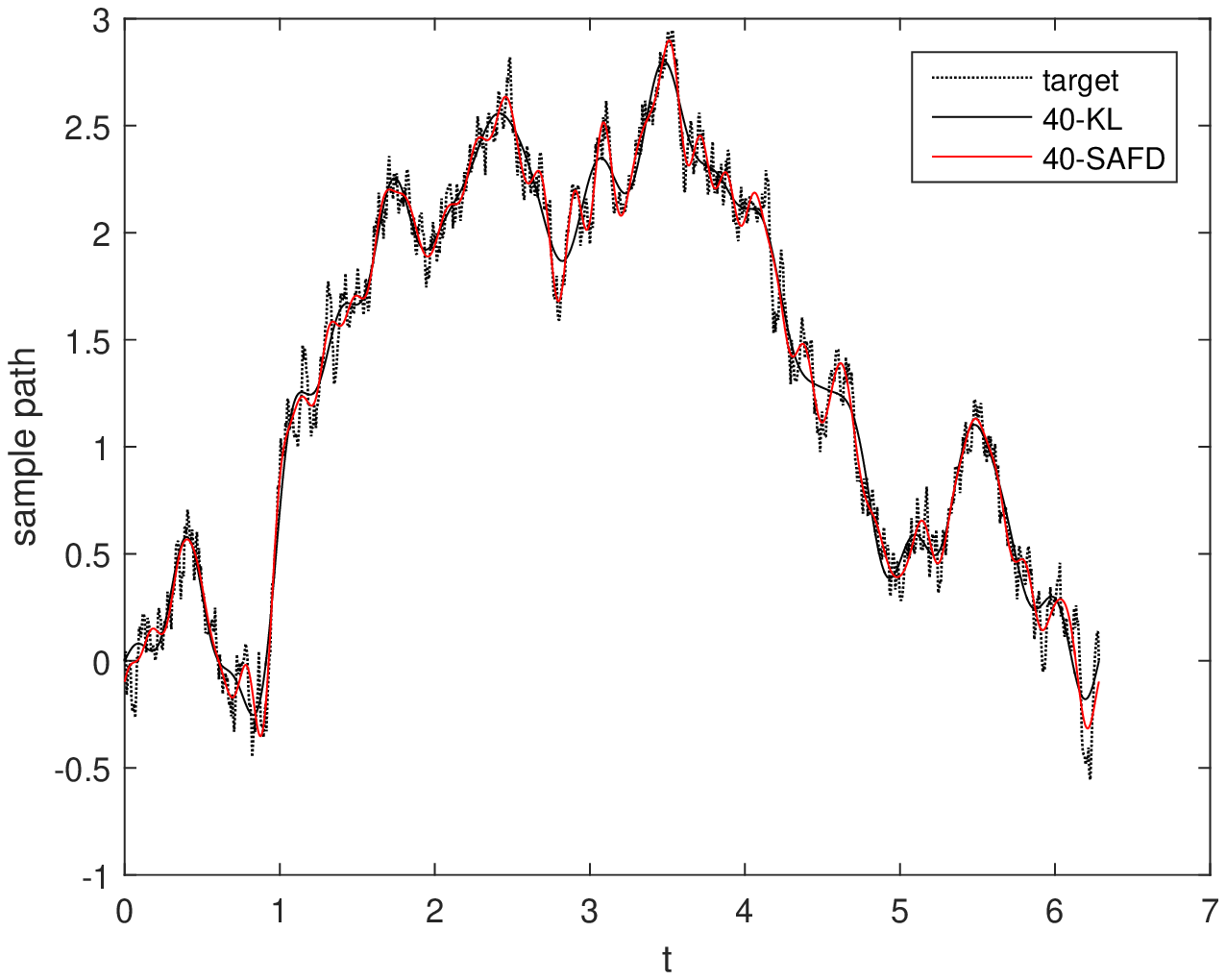}
		\centerline{{\scriptsize SAFD: 40 partial sum}}
	\end{minipage}
	\caption*{Figure 2b: sample path IIb ($2^{10}$ points) of Brownian bridge}
\label{figureB2b}
\end{figure}
\begin{table}[h]
\scriptsize
  \centering
    \begin{tabular}{c|c|c|c|c}
    \toprule
    \multicolumn{1}{c|}{$n$ partial sum} & \multicolumn{1}{c|}{$n=10$} & \multicolumn{1}{c|}{$n=20$} & \multicolumn{1}{c|}{$n=30$} & \multicolumn{1}{c}{$n=40$} \\
    \midrule
    KL    & 0.0245 & 0.0103 & 0.0074 & $0.0059$ \\
    \midrule
    SAFD & 0.0120 & 0.0061 & 0.0046 & $0.0035$ \\
    \bottomrule
    \end{tabular}%
    \caption*{Table 2b: Relative error}
      \label{tableSAFD2}%
\end{table}%

\begin{example}(SnB on the Szeg\"o kernels dictionary)\label{Brownian3}
The Brownian bridge is generated by using 2048 sampling points in $[0,2\pi]$ with the uniform spacing $\Delta t \thickapprox 0.003$. In this example, we approximate the sample path with the KL expansion and the SnB method. As shown in Figure 3 and Table 1, with all the 15, 30, 60, 100 partial sum approximations the SnB method outperforms the KL method in the details. The relative errors are given in Table 3.
\end{example}

\begin{figure}[H]
	\begin{minipage}[c]{0.23\textwidth}
		\centering
		\includegraphics[height=4cm,width=3.3cm]{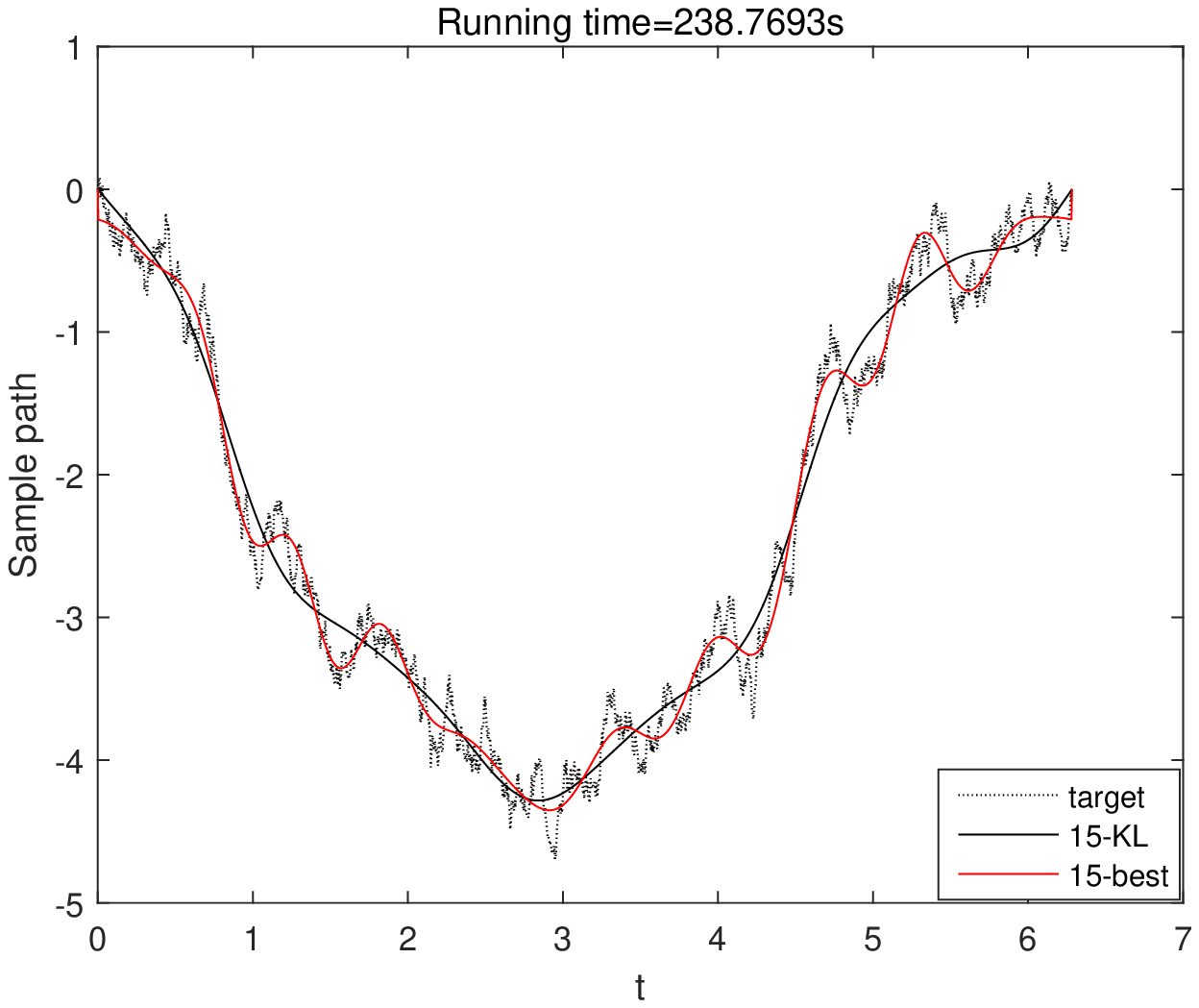}
		\centerline{{\scriptsize 15-best}}
	\end{minipage}
	\begin{minipage}[c]{0.23\textwidth}
		\centering
		\includegraphics[height=4cm,width=3.3cm]{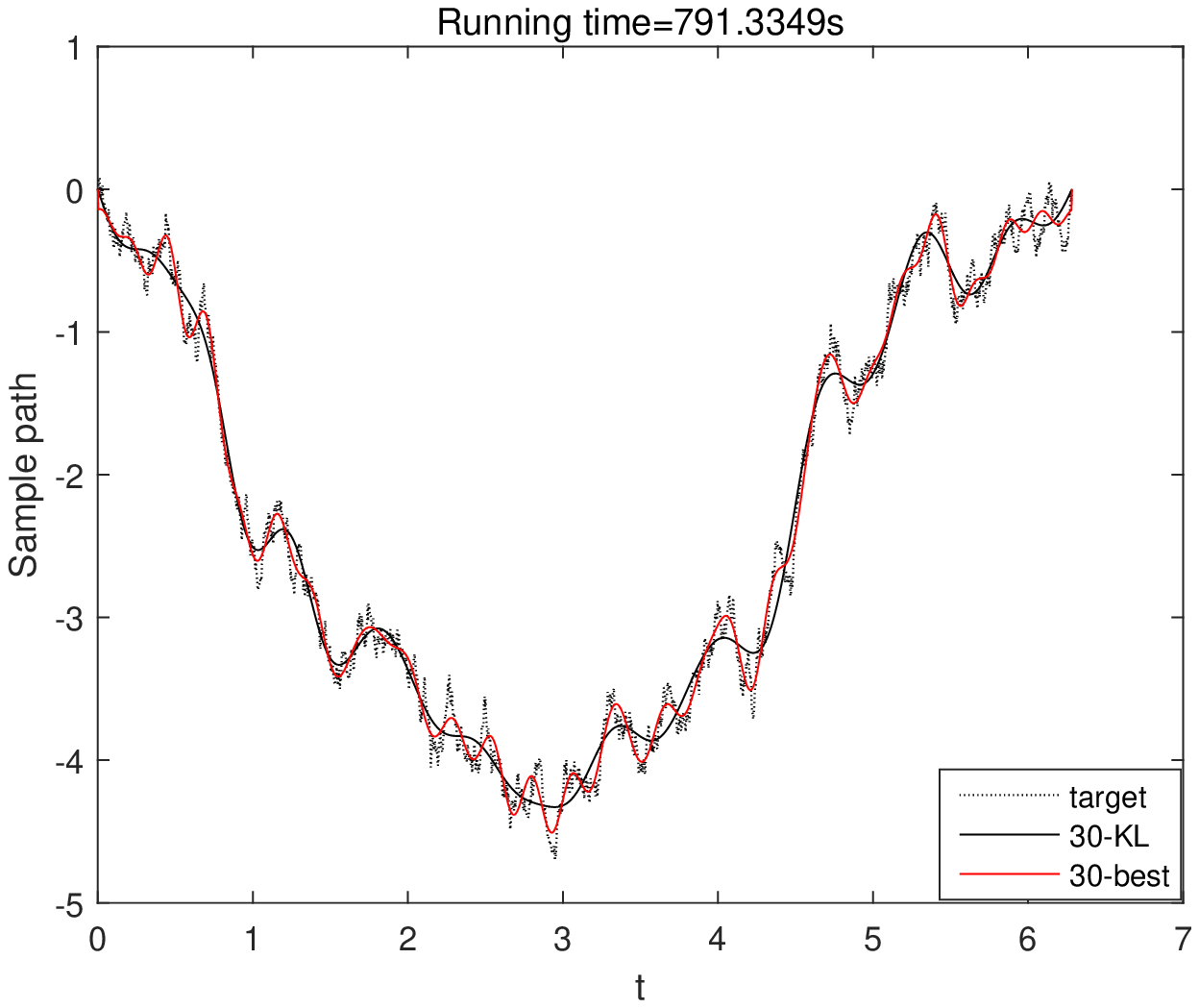}
		\centerline{{\scriptsize 30-best}}
	\end{minipage}
	\begin{minipage}[c]{0.23\textwidth}
		\centering
		\includegraphics[height=4cm,width=3.3cm]{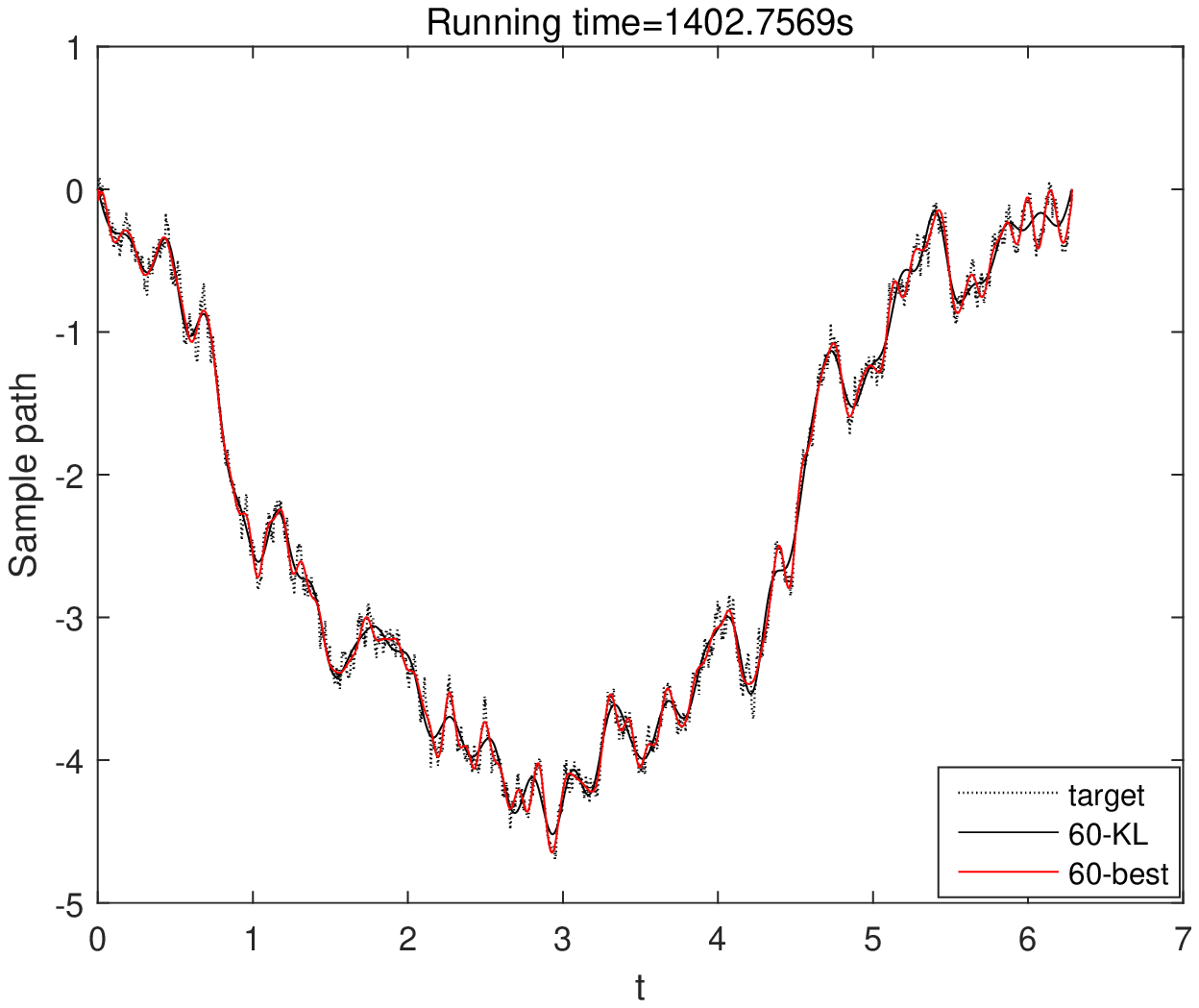}
		\centerline{{\scriptsize 60-best}}
	\end{minipage}
	\begin{minipage}[c]{0.23\textwidth}
		\centering
		\includegraphics[height=4cm,width=3.3cm]{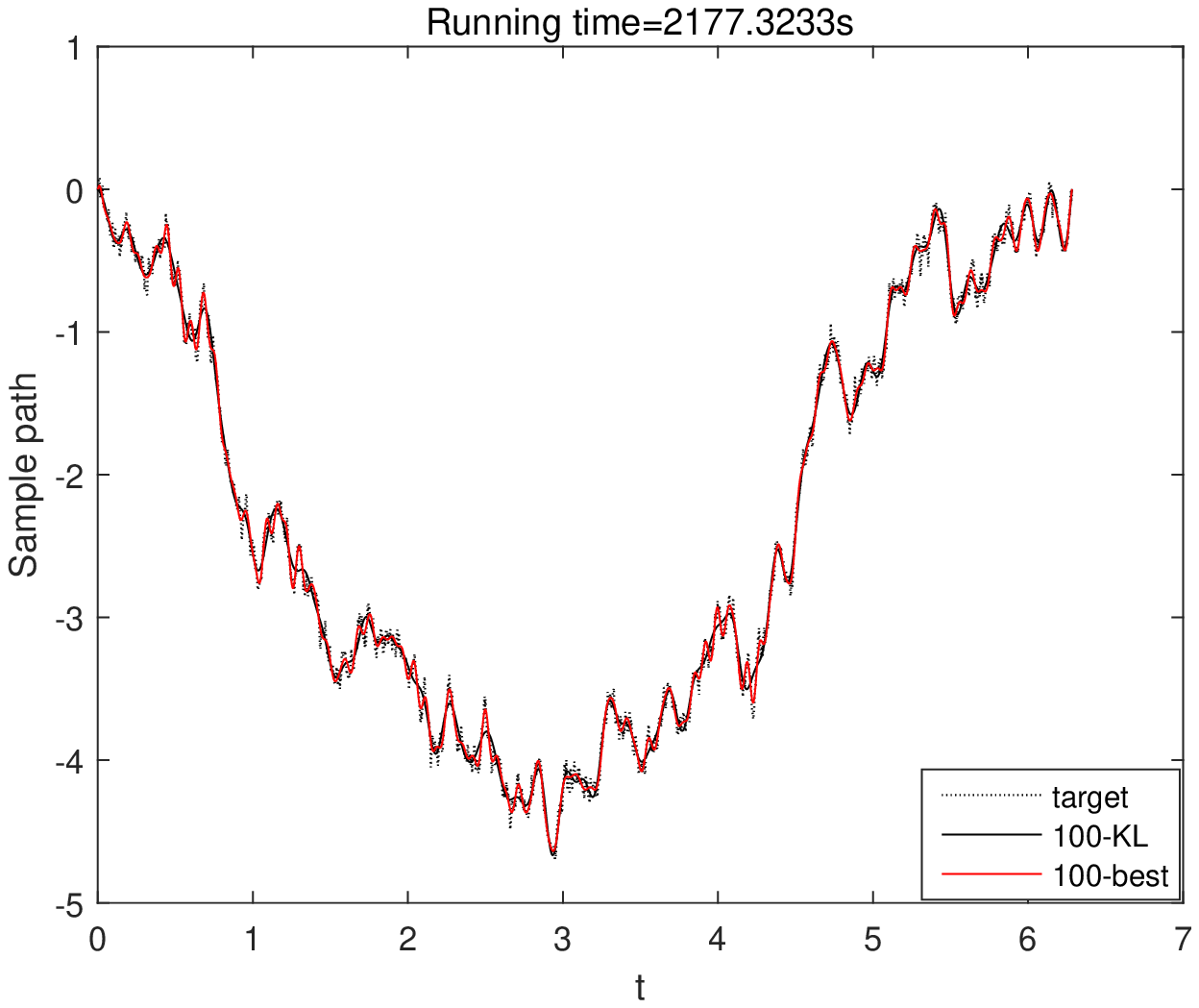}
		\centerline{{\scriptsize 100-best}}
	\end{minipage}
	\caption*{Figure 3: sample path III ($2^{11}$ points) of Brownian bridge}
\label{figureB3}
\end{figure}
\begin{table}[H]%
\tiny
  \centering
    \begin{tabular}{c|c|c|c|c}
    \toprule
           \multicolumn{1}{c|}{$n$ partial sum} & \multicolumn{1}{c|}{$n=15$} & \multicolumn{1}{c|}{$n=30$} & \multicolumn{1}{c|}{$n=60$} & \multicolumn{1}{c}{$n=100$} \\
    \midrule
    KL    & 0.0068 & 0.0031 & 0.0015 & 8.2008$\times10^{-4}$ \\
    \midrule
    SnB & 0.0031 & 0.0015 & 6.9197$\times10^{-4}$ & 3.8540$\times10^{-4}$ \\
    \bottomrule
    \end{tabular}%
    \caption*{Table 3: relative error}
      \label{table1}
\end{table}%

\section{Conclusions}
In the article we establish, by using the covariant function, the algorithms of SAFD, SPOAFD and SnB, and prove that they enjoy the same convergence rate as it does by the KL decomposition method. The experimental examples show that the AFD type methods perform even better than the KL method in describing the local details.  The AFD type methods gain extra adaptivity from dictionary selection according to the tasks taken. To perform the algorithms of SAFD, SPOAFD and SnB one does not need to compute eigenvalues and the eigenfunctions of the integral operator defined by the covariance kernel that, compared with the KL method, greatly reduces computation complexity and save computer consumes.  The proposed AFD type expansions are well applicable also to infinite time intervals, as well as to unbounded regions for the space variable.

\section*{Acknowledgement}
The authors wish to sincerely thank Professor Radu Victor Balan whose relevant questions inspired this study.

\end{document}